\documentclass{article}
\usepackage[utf8]{inputenc}
\usepackage[centertags]{amsmath}
\usepackage{caption,graphicx}
\usepackage{subcaption}
\usepackage[text={6.0in,8.6in},centering,letterpaper]{geometry}
\usepackage[numbers,comma,square,sort&compress]{natbib}
\usepackage[utf8]{inputenc}
\usepackage{mathrsfs}
\usepackage{amsthm, amsfonts, mathrsfs, amsfonts, amssymb}
\usepackage[british]{babel}
\usepackage{mathtools}
\usepackage{comment}
\usepackage{url} 
\usepackage{xcolor}

\usepackage{authblk}
\usepackage{soul}
\usepackage{todonotes}

\everymath={\displaystyle}
\numberwithin{equation}{section}

\newcommand{\R}{\ensuremath{\mathbb{R}}}

\newcommand{\p}{\ensuremath{\partial}}
\newcommand{\s}{\ensuremath{\sigma}}
\newcommand{\tr}{\mathrm{tr}} 

\newcommand{\sref}[1]{(\ref{#1})}                       

\newtheorem{remark}{Remark}
\newtheorem*{defn*}{Definition}

\title{Blurring the Busse balloon: Patterns in a stochastic Klausmeier model.}

\begin{document}
\author[1,2,3]{Christian H.S. Hamster}
\affil[1]{Dutch Institute for Emergent Phenomena, University of Amsterdam, Amsterdam, The Netherlands}
\affil[2]{Korteweg-De Vries Institute for Mathematics, University of Amsterdam, Amsterdam, The Netherlands}
\affil[3]{Biometris, Wageningen University \& Research, Wageningen, The Netherlands}

\author[3]{Peter van Heijster}

\author[3]{Eric Siero}

\maketitle

\begin{abstract}
We investigate the effect of stochastic forcing on the stable periodic patterns of the Busse balloon in a one-dimensional Klausmeier model for dryland vegetation patterns. Using numerical methods, we can accurately describe the transient dynamics of the stochastic solutions and compare several notions of stability. In particular, we show that 
the boundary of the Busse balloon becomes blurred under the stochastic perturbations and that the
stochastic stability heavily depends on the model parameters, the intensity of the noise and the location of the wave number of the periodic pattern within the deterministic Busse balloon.   \\
\end{abstract}

\section{Introduction}
\begin{figure}[t]
  \centering
 		\def\svgwidth{0.7\columnwidth}
    		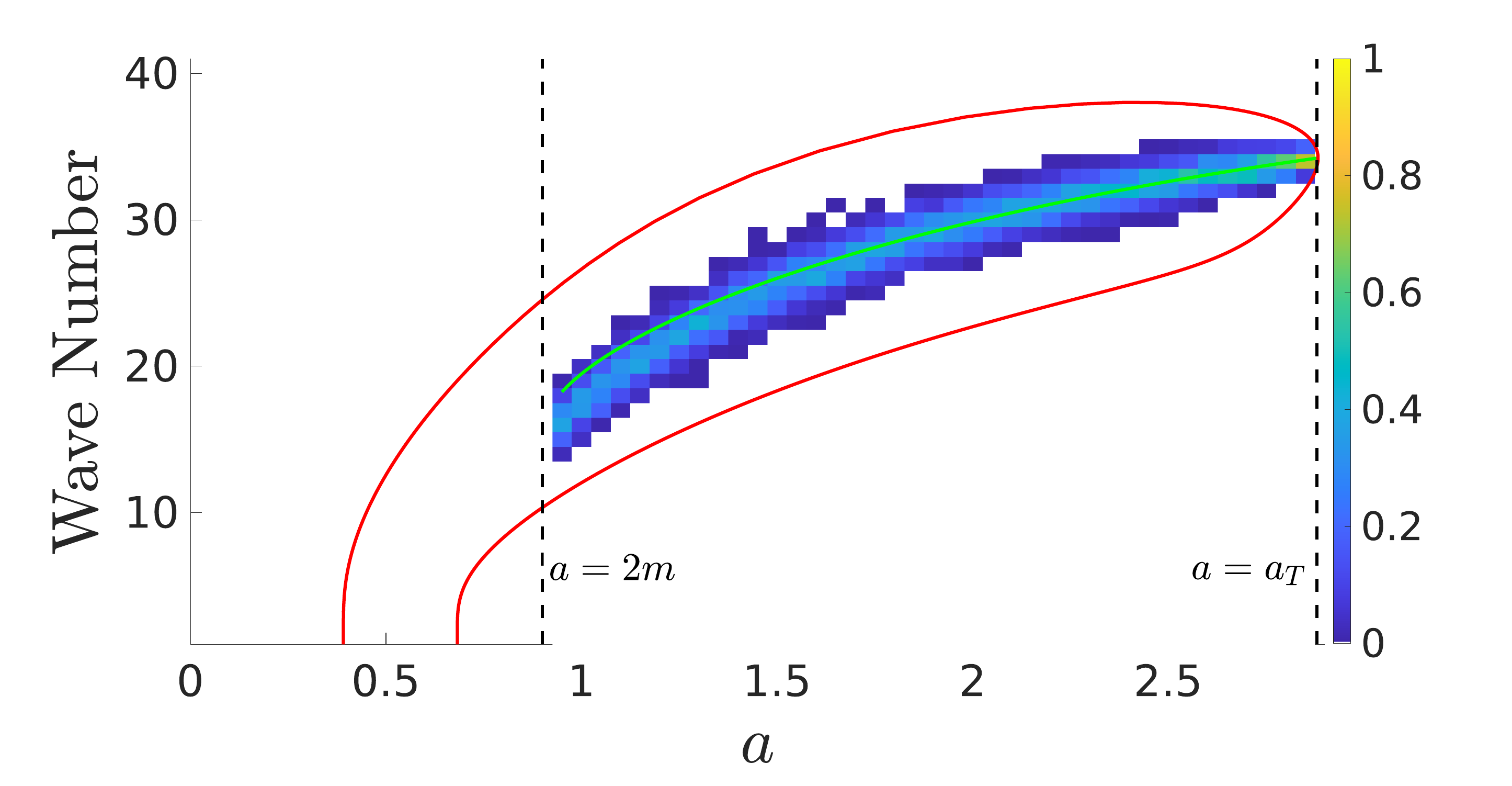
\caption{We show the normalized histogram of $200$ runs per $a-$value of the deterministic version of \eqref{eq:int:StochK} with $\sigma =0$, $m=0.45$ and $d=500$ on the periodic domain $[-250,250]$. The runs are initiated from the Turing unstable spatially homogeneous steady state \eqref{ss}, which exists for $a\geq 2m=0.9$, plus a small perturbation. The vertical axis indicates the predominant wave number at the end of each simulation ($T=5000$). Wave numbers that are not observed are shown as white. The area enclosed by the red curve (computed with \texttt{auto-07p}~\cite{auto}) corresponds to the Busse balloon of all stable periodic patterns. The green line indicates the most unstable wave number, see Appendix \ref{app:MostUnst}.}
\label{fig:int:FromUnif}
\end{figure}

Self-organization leads to spatially periodic patterns in nature. Examples of patterns in ecology can be found in drylands~\cite{Macfadyen1950,deblauwe2008global, bastiaansen2018multistability}, mussel beds~\cite{vandeKoppel2008experimental} and peatlands~\cite{Fosteretal1983,couwenberg2005self}. These pattern formation phenomena are often modelled with Reaction-Diffusion Equations (RDE) and they have been studied extensively, mainly in the deterministic setting~\cite[e.g.]{Murray2003MathematicalBiologyII,kondo2010reaction}. The deterministic RDE models often support multiple stable and unstable spatially-periodic stationary patterns for a given set of parameters. The stable periodic patterns of these models are graphically represented via their wave numbers in so-called Busse balloons (with one of the system parameters on the horizontal axis, see Figure~\ref{fig:int:FromUnif} for an example). Busse balloons were originally introduced in the field of fluid mechanics to represent stable patterns of convection rolls~\cite{busse1978non}. However, nowadays, they are generic objects for spatially patterned systems.

The Busse balloon can thus be used to predict which patterns are observable in a system. However, in any real-world system, noise will perturb these patterns. If the noise intensity is small, stochastic solutions will generally remain in the basin of attraction of deterministically stable patterns on long timescales, so the deterministic Busse balloon suffices. But for larger noise intensities this is not the case and the notion of stability, observability, and the Busse balloon needs to be re-evaluated.

In this article, we propose a framework to address the above-raised questions for RDE models influenced by noise, {\emph{i.e.}} for stochastic RDE models. Although the techniques in this article are applicable to any stochastic RDE model that possesses a Busse balloon (in the deterministic setting), we focus on a model that models pattern formation in drylands. Drylands make up 41\% of the terrestrial land surface, with a human population of over two billion people~\cite{MillenniumEcosystemAssessment}. The presence of spatial vegetation patterns in these areas is ubiquitous~\cite{deblauwe2008global}. The patterns have on the one hand been interpreted as early warning signals for critical transition~\cite{Schefferetal2009} and on the other hand have been linked to ecosystem resilience~\cite{Rietkerketal2021}, see also the recent review in~\cite{martinez2023integrating}. In \cite{KASTNER2024}, the regularity of dryland vegetation patterns is discussed and it is shown that patterns in deterministic RDEs are too regular when compared with real vegetation patterns, which implies that stochastic modelling is necessary.  

Since the creation of the first multi-component deterministic model for vegetation patterns by Klausmeier in '99~\cite{klausmeier1999}, which forms the basis of our stochastic model \eqref{eq:int:StochK} below, many extensions of this model have been proposed. The equations can be modelled in two space dimensions and with a (possibly varying) hill slope~\cite{bastiaansen2019stable}, to include grazing~\cite{siero2018grazing} or fluctuations in rainfall~\cite{eigentler2020effects}, to name a few. Similar models with three instead of two components include~\cite{HilleRisLambersetal2001,Giladetal2004}. 

We study the influence of stochastic perturbations, or noise, on the stability and observability of periodic patterns in a one-dimensional non-dimensionalized Klausmeier model for dryland vegetation patterns on a flat terrain~\cite{klausmeier1999,siteuretal2014}:
\begin{align}
    \begin{split}
        du =& [du_{xx}+a-u-uv^2]dt\,,\\
        dv =& [v_{xx}-mv+uv^2]dt+\sigma vdW^Q_t\,.
    \end{split}
    \label{eq:int:StochK}
\end{align}
The equation is posed on a spatial domain $[-L,L]$ with periodic boundary conditions for $t\in\R^+$.
The (nondimensionalized) biomass of the system is represented by $v(x,t)$ and $u(x,t)$ represents the (nondimensionalized) amount of water. The diffusion of water (in relation to the diffusion of the biomass) is modelled by the diffusion coefficient $d$ and is thus assumed to be large, $a$ models the (average) rainfall and $m$ is the (average) mortality rate of the biomass.
Here, we assume that this mortality rate is noisy in nature, for example, due to human disturbances or grazing (however, see Remark~\ref{R:noise}). That is, we replace the mortality rate $m$ of the deterministic model by $m+\sigma dW^Q_t$, which results in the multiplicative noise term (in It\^o interpretation) $\sigma v dW^Q_t$ in \eqref{eq:int:StochK}. Here, $W^Q_t$ is a Gaussian process on $L^2([-L,L])$ with no correlation in time, {\emph{i.e.}} white in time, and the spatial correlation is given by a Gaussian function $q(x-y)$, {\emph{i.e.}} coloured in space. This means that correlations of increments of the stochastic process $W^Q_t$ are given by $E[dW_t(x)dW_s(y)]=\delta(t-s)q(x-y)$, with $\delta(\cdot)$ the Dirac-delta function. See~\cite{lord2014book} for an introduction to stochastic RDEs and related numerical techniques.

\begin{remark}
\label{R:noise}
In \eqref{eq:int:StochK} we assumed the mortality rate $m$ of the biomass to be the only noisy term in the model. We could also have included other stochastic processes in the model. For example, we could replace the rainfall $a$ with a stochastic process that actually models rainfall patterns in, say, Somalia~\cite{gandhi2023pulsed}. However, as the Klausmeier model is a conceptual model, including such detail in the model would force us to also significantly extend the other terms in the model (as was done in~\cite{gandhi2020fast,gandhi2023pulsed}). 
\end{remark}

\paragraph{Deterministic RDE model: stability and observability}
The deterministic version of \eqref{eq:int:StochK}, {\emph{i.e.}} \eqref{eq:int:StochK} with $\sigma=0$, has up to three stationary spatially homogeneous solutions: the bare-soil state $(u,v)=(a,0)$ exists for all $a$, and, for $a\geq 2m$, two more background states appear from a saddle-node bifurcation. From these, the state $(u,v)=(a/2+\sqrt{a^2-4m^2}/2,(a-u)/m)$ is unstable and 
\begin{align}
\label{ss}(\bar u,\bar v)=(a/2-\sqrt{a^2-4m^2}/2,(a-\bar u)/m)
\end{align} is the one that is subjected to a supercritical Turing instability against a critical wave number $k=k_T$ at some $a=a_T$. For $a\leq a_T$, a family of stable and unstable periodic patterns can be found, see \cite{vdStelt2013rise} for a discussion on several ways in which a pattern can destabilize. Starting at the Turing bifurcation, the boundary between stable and unstable patterns can be traced through parameter space by means of continuation~\cite{RademacherSandstedeScheel2007}. This way, the region in parameter space of all stable periodic patterns can be determined and this region is known as the Busse balloon. For the deterministic version of \eqref{eq:int:StochK}, this results in a boundary given by the red curve in Figure~\ref{fig:int:FromUnif}.

Given that for a given set of parameters, multiple patterns are stable, it is natural to ask if all patterns are {\emph{equally likely}} to be observed, either in real ecosystems or in numerical simulations of the deterministic RDE model. 
This question of observability highly depends on the initial condition of the RDE simulation and because of the many alternative stable states, the likeliness to observe a specific pattern is history-dependent~\cite{sherratt2013history}. 
When we start close to a periodic pattern with stable wave number $n$, we expect to observe the same wave number at the end of the simulation. However, when we start at an unstable homogeneous background state \eqref{ss} plus a perturbation, it is hard to predict which wave number will be observed at the end of the simulation.
As the initial condition is random, we can expect different results in each simulation. Indeed, as Figure~\ref{fig:int:FromUnif} shows, different wave numbers (see Remark~\ref{R:WN}) are selected, but the distribution is far from uniform within the Busse balloon (represented by the red curve in the figure). A large part of the balloon is not selected, or selected with such a small probability that we do not observe it in a simulation with 200 runs. We observe a peak around the \emph{most unstable mode} (green curve in Figure~\ref{fig:int:FromUnif}), as was observed before~\cite{sherrattlord} in a Klausmeier model with advection. See Appendix~\ref{app:MostUnst} for a computation of the most unstable mode.

\begin{remark}
\label{R:WN}
For a periodic pattern with wavelength $\lambda$, the wave number is defined as $1/\lambda$. Given a periodic domain $[-L,L]$, a periodic pattern has $2L/\lambda$ maxima or, in the ecological interpretation, $2L/\lambda$ patches of vegetation. Throughout this whole article, we always choose $L=250$, hence without risk of confusion, we use the term `wave number' to mean $2L/\lambda$. This allows us to directly talk about patterns having, say, a wave number of $29$ or $30$, instead of 0.058 or 0.06. Furthermore, we can now directly compare the wave number to the pulse number, {\emph{i.e.}} the number of pulses we find by counting the number of maxima. 
\end{remark}

\paragraph{Stochastic RDE model: first exit time and local wave number}
The stationary spatially-periodic patterns of the deterministic RDE model are no longer stationary solutions of the stochastic RDE. Furthermore, due to the structure of the multiplicative noise, only the homogeneous bare-soil state $(u,v) = (a,0)$ remains as background state. In a pessimistic view, one might say that this makes the notion of a {\emph{stochastic Busse balloon}} void. 
However, as the numerical simulations will show, initial conditions that start close to a deterministic pattern will remain close to this pattern up to a certain random time, which we refer to as the {\emph{first exit time}}. 
When a deterministically unstable pattern is chosen as initial condition, this first exit time can be very short, but for a deterministically stable pattern, it can be long, especially for small $a$- and $\sigma$-values. In other words, the concept of the first exit time in the stochastic case mimics the concept of stability from the deterministic case. Therefore, we use the {\emph{average first exit time}} to study the stability for periodic-like patterns in stochastic RDEs. Note that the concept of average first exit times is, for example, also used in~\cite{arani2021exit} as a measure of ecological resilience.

Average first exit times give an indication of the stability of the patterns but do not say anything about what happens after destabilization (similar to the deterministic case: the Busse balloon is related to stability, but not necessary observability, see Figure~\ref{fig:int:FromUnif}). If the solution after destabilization is close to a periodic pattern, it will again destabilize after a certain waiting time. However, there is no guarantee that the solution is actually close to any of the periodic solutions. Hence, in order to talk about observability, {\emph{i.e.}} what are the typical properties of solutions we observe after a certain integration time, we must also be able to classify solutions that are not close to periodic patterns. 

A first natural candidate for this is the {\emph{predominant wave number}}, {\emph{i.e.}} the wave number where the absolute value of the Fourier transform attains its maximum. However, the dynamics from one state, say with $n$ pulses, to another state, say with $n-1$ pulses, is -- already in the deterministic case -- slow, see Figure~\ref{fig:DetLocal} for an example. Especially, the predominant wave number will only change from $n$ to $n-1$ after a long integration time. Therefore, it is {\emph{a priori}} not clear what wave number should be assigned to the solution during this {\emph{transient phase}}. In the limit of $t\to\infty$, this transient dynamics is not a problem, but on the short timescales of the stochastic dynamics, it is. Therefore, the questions about observability will depend on what definition of wave number we assign to a transient state. 

Another natural candidate for classification would be to count the number of maxima, which we refer to as the pulse number. However, just the pulse number does not tell us how these pulses are distributed over the domain, or how close the patterns are to being periodic. In Figure~\ref{fig:DetLocal}, the pulse number is fixed at $29$, but the dynamics is driven by the fact that, initially, the pulses are not spread out evenly. For a discussion of the dynamics of unevenly spread out pulses in the deterministic setting, see \cite{bastiaansen2019dynamics}.

In the context of fluid dynamics, where the Busse balloon originates, the concept of a local wavelength through a local Fourier transform is used to describe stochastic pattern formation in a Swift-Hohenberg equation~\cite{vinals1991numerical}. We will also use this concept of {\emph{local wave numbers}} to study the observability of periodic-like patterns in stochastic RDEs, as we will show that the local wavenumber captures both the quasi-periodicity of the patterns and the number of pulses.

This article is organized as follows. In \S\ref{sec:NumProp}, we further introduce and define the stochastic concepts of (average) first exit time and local wave number to study the stability and observability of periodic-like patterns in stochastic RDEs. Subsequently, in \S\ref{results}, we use these concepts to blur the Busse balloon: we show that the average first exit time of a periodic-like pattern of \eqref{eq:int:StochK} heavily depends on the position of that pattern within the Busse balloon and we use local wave numbers to describe the stationary distribution of the stochastic dynamics. We end the article in \S\ref{sec:disc} with a discussion of the results and discuss possible avenues of further research.

\section{Methods}
\label{sec:NumProp}
We will study the solutions of the stochastic RDE \sref{eq:int:StochK} numerically. Especially, we will focus on classifying solutions using several methods, such as the average first exit time and local wave numbers, and compare their similarities and differences. 
Using the continuation software \texttt{pde2path}~\cite{uecker2014pde2path,uecker2022continuation} we first compute for the deterministic case a set of stable and unstable periodic patterns for a discrete set of rainfall values $a$ (multiples of $0.05$) and wave numbers. We use these patterns as initial conditions for the stochastic simulations.
The numerical scheme behind these stochastic simulations is based on~\cite{lord2014book} and can be found in Appendix~\ref{app:NumScheme}. The remaining system parameters stay fixed throughout this article (unless mentioned otherwise) and their values can be found in Table~\ref{tab:Par}.

\begin{table}
    \centering
    \begin{tabular}{|c|c|c|c|}
    \hline
    Parameter & Meaning & Value & Reference\\
    \hline \hline
        $d$ & Diffusion coefficient of water & 500 &~\cite{siteuretal2014} \\
         $m$& Average mortality of vegetation  & 0.45&~\cite{klausmeier1999}\\
         $L$& Length of domain  & 250 & -\\
         $\ell$& Window width in local Fourier transform& 50& -\\
         \hline
    \end{tabular}
    \caption{Table of fixed system parameters.}
    \label{tab:Par}
\end{table}

First, we will show some typical realizations of the stochastic solutions and then discuss the question of classification, {\emph{i.e.}} how do we determine the number of pulses in the solution and to what extent is this related to the predominant wave number. Of course, `typical' depends on the chosen parameters, especially $\sigma$. See Fig.~\ref{fig:app:Sigma} for a comparison of different noise intensities.  

\subsection{Single realizations and their wave numbers}
\begin{figure}[b]
\begin{subfigure}{.49\textwidth}
  \centering
 		\def\svgwidth{\columnwidth}
    		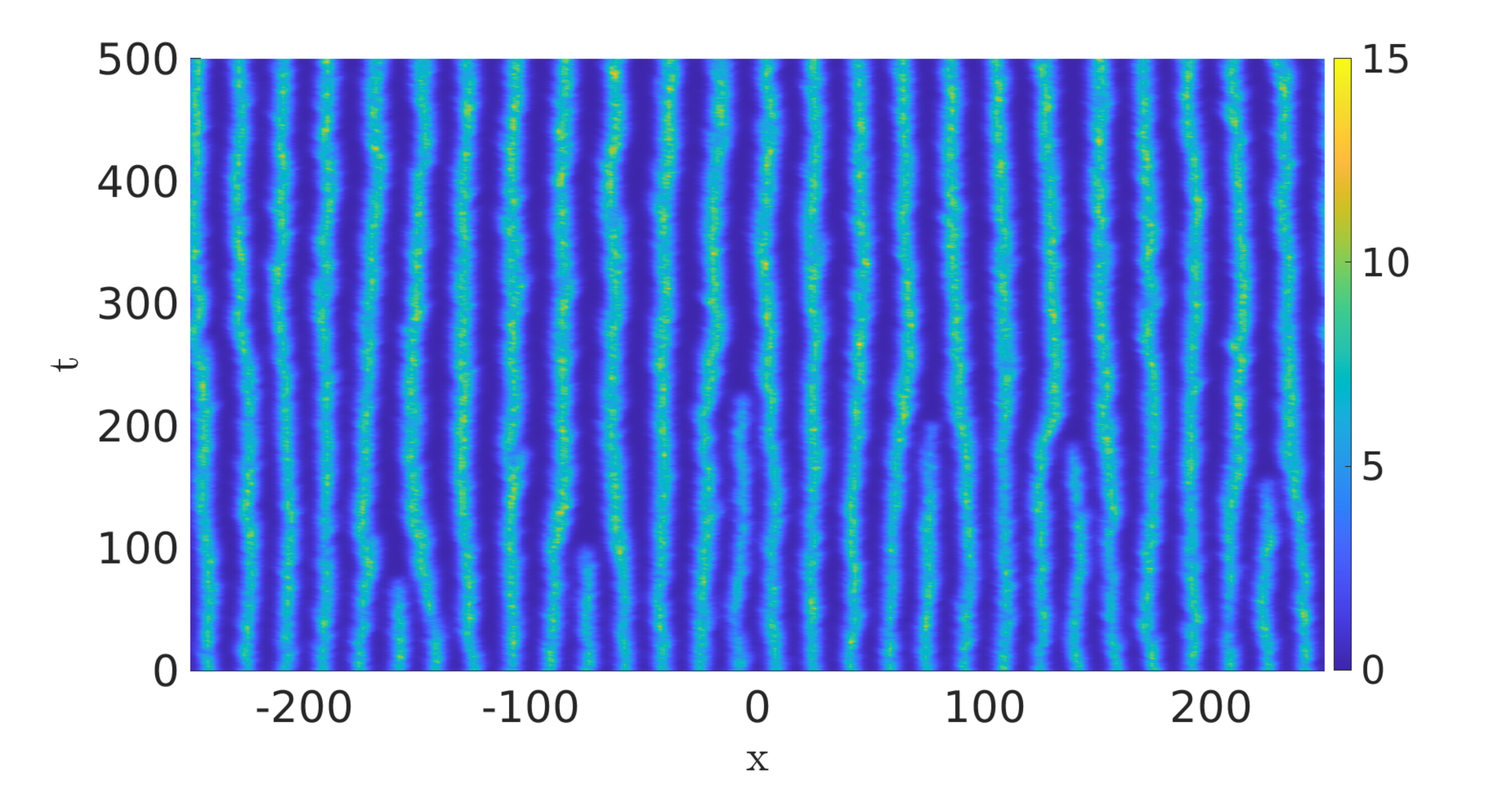
  \caption{}
    \label{fig:ExampleVn30a15}
\end{subfigure}
\begin{subfigure}{.49\textwidth}
  \centering
 		\def\svgwidth{\columnwidth}
    		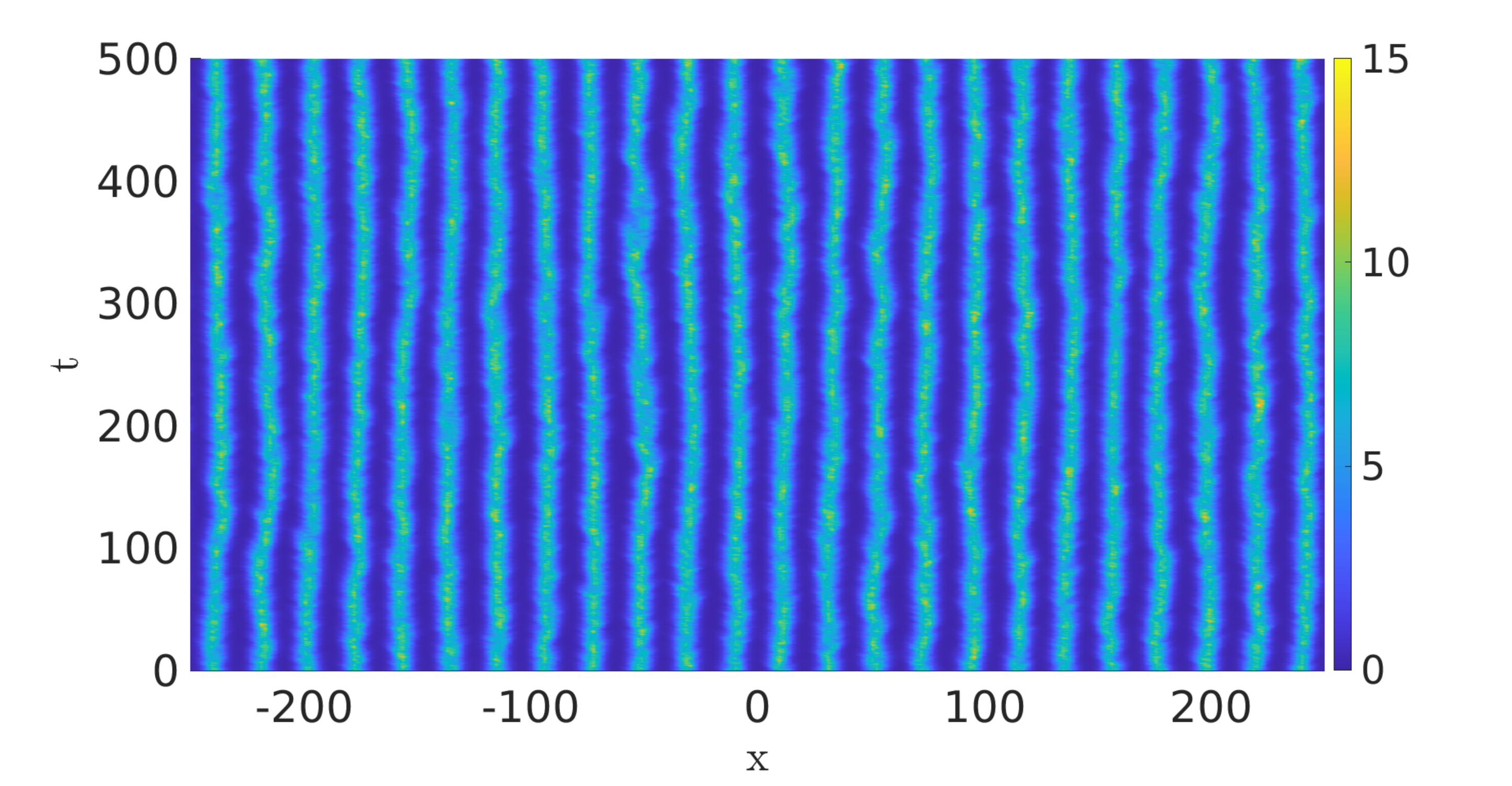
  \caption{}
    \label{fig:ExampleVn24a15}
\end{subfigure}
\caption{Single realization of the stochastic RDE~\sref{eq:int:StochK} for $a=1.5$ and $\sigma=0.2$, starting with a periodic solution with (a)~$30$ pulses and (b)~24 pulses. Note that only the $v$-component ({\emph{i.e.}} vegetation) is shown. In (a), six out of the $30$ pulses disappear due to the stochasticity, while none disappear in (b). In (a), the first exit time $T_\mathrm{exit}\approx 74$, while in (b) $T_\mathrm{exit} = T_{\rm max}=500$.}
\label{fig:SingleR}
\end{figure}

In Figure~\ref{fig:SingleR}, we show two realizations of the stochastic RDE~\sref{eq:int:StochK} with $a=1.5$ and $\sigma=0.2$ (and with the remaining system parameters as given in Table~\ref{tab:Par}). In Figure~\ref{fig:ExampleVn30a15}, the initial condition is the stable periodic solution (in the deterministic case) with wave number $30$, whilst the initial wave number of the initial condition is $24$ in Figure~\ref{fig:ExampleVn24a15}. Even though the deterministic solution with wave number~$30$ is stable, we observe that several pulses disappear due to the stochastic forcing, resulting in a pattern with $24$ pulses. Does this imply that wave number $24$ is in some sense more stable than wave number $30$? Indeed, starting from wave number $24$ in Figure~\ref{fig:ExampleVn24a15} we do not observe a change in the number of extrema. However, this might change from realization to realization and only repeated numerical experiments can tell us if the dynamics observed in Figure~\ref{fig:SingleR} is in some sense typical.

\begin{figure}[t]
\begin{subfigure}{.49\textwidth}
  \centering
 		\def\svgwidth{\columnwidth}
    		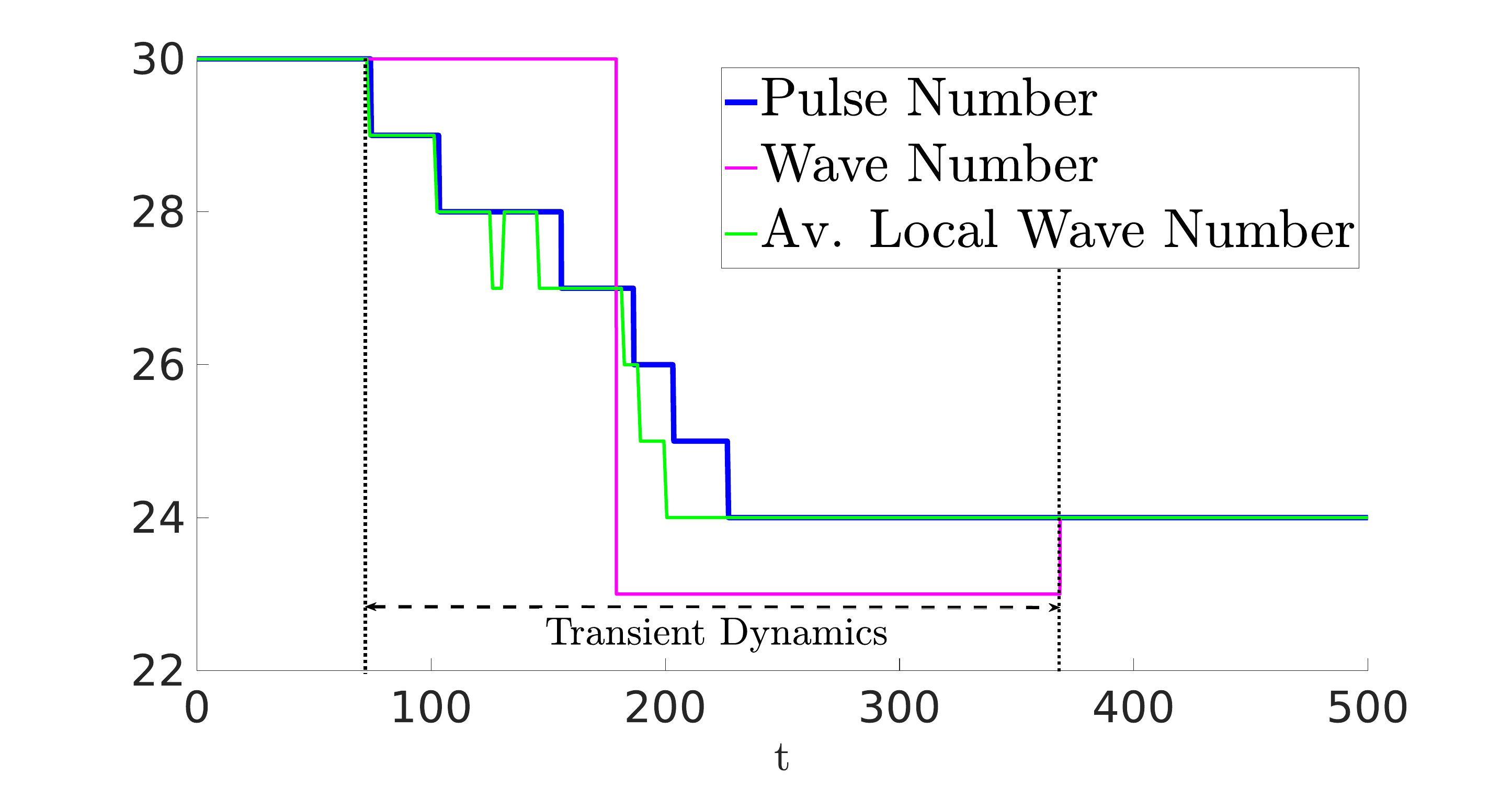
  \caption{}
    \label{fig:ExampleNumbern30a15}
\end{subfigure}
\begin{subfigure}{.49\textwidth}
  \centering
 		\def\svgwidth{\columnwidth}
    		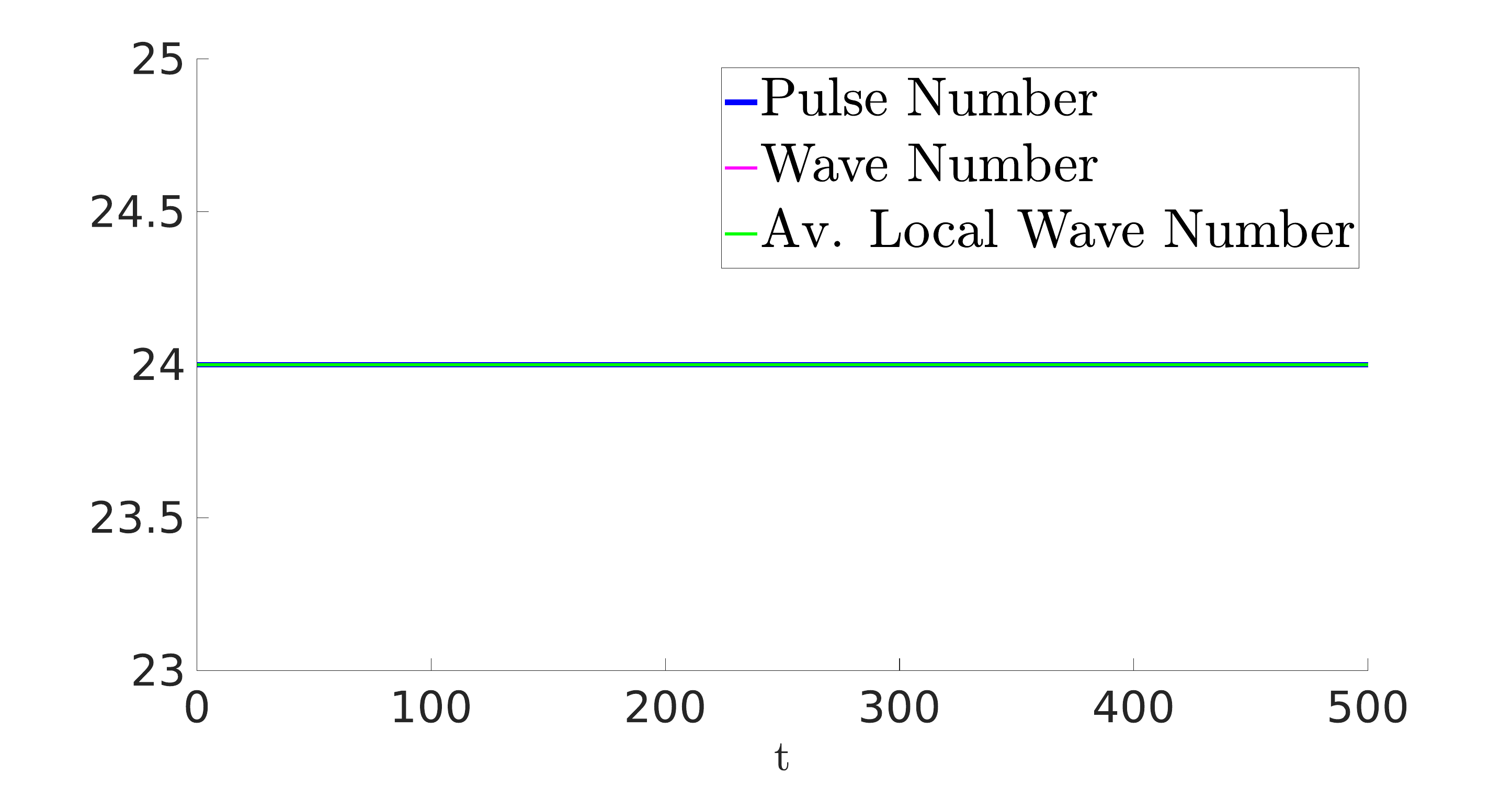
  \caption{}
    \label{fig:ExampleNumbern24a15}
\end{subfigure}
\caption{Figure (a) shows the predominant wave number, the pulse number, and the rounded average local wave number corresponding to Figure~\ref{fig:ExampleVn30a15}, while~(b) shows these three numbers for Figure~\ref{fig:ExampleVn24a15}. In~(b), the curves overlap.}
\label{fig:WNvsPN}
\end{figure}

For a stationary spatially-periodic pattern, the number of pulses and predominant wave number are equal. However, when we delete one pulse in the pattern, these two numbers are not the same anymore. The predominant wave number of the pattern has not changed, but the number of pulses clearly did decrease by one. When we let time evolve and the solution converges to a (different) stationary spatially-periodic solution, then the pulse number and predominant wave number coincide again. Hence, there is no discrepancy between these two measures for stationary solutions, but in the {\emph{transient phase}}, the two numbers can differ. More precisely, we define the transient phase of a solution as the time interval where the predominant wave number and pulse number differ.
In the stochastic version, the periodic patterns are not stationary solutions anymore, hence there can always be a discrepancy between the wave number and pulse number. Figure~\ref{fig:WNvsPN} shows the predominant wave number and number of pulses associated with the two realizations in Figure~\ref{fig:SingleR}. In Figure~\ref{fig:ExampleNumbern24a15} these numbers coincide. However, in Figure~\ref{fig:ExampleNumbern30a15}, the predominant wave number makes one large step from $30$ to $23$ and then steps up again to $24$ while the pulse number decreases in single steps to $24$. Hence, the classical Fourier transform used to compute the predominant wave number has its limitations in describing transient states. See Appendix~\ref{app:numbers} for an explanation of the numerical classification of the pulse number and wave number.

\subsection{Local wave numbers}
In studies on dryland patterns using $2$D remote sensing data such as~\cite{bastiaansen2018multistability}, the total observed area with spatial patterns is divided into smaller squares, and in each square, the predominant wave number is computed. Similarly,  
we can apply a window to the domain $[-L,L]$ to study the local properties of the solution to further analyze the transient states and address the observability of the patterns. 

Following~\cite{vinals1991numerical}, we define a local Fourier series 
\begin{align}
    \mathcal{F}_{x_0}[u](k,t)=\frac{1}{2L}\int_{-L}^Lu(x,t)e^{\frac{(x-x_0)^2}{\ell^2}}e^{\frac{2\pi ikx}{L}}dx,
\end{align}
where $\ell$ defines the width of the window. The predominant mode, {\emph{i.e.}} the value of $k$ for which $|\mathcal{F}_{x_0}[u](k,t)|^2$ attains its maximum, now depends on $\ell$, $t$ and $x_0$. 
The information in the local wave number can be used in two ways to obtain information per timestep. First, we can use the local wave numbers to plot (per timestep) a histogram, that indicates the distribution of the local wave numbers. This allows us to characterize a state based on the shape of its distribution. Secondly, we can use the distribution to assign a single number to the solution at each point in time, {\emph{e.g.}} we can compute the average local wave number rounded to the nearest integer. This gives us a notion of a {\emph{typical local wave number}}. As the green line in Figures \ref{fig:ExampleNumbern30a15} and \ref{fig:ExampleNumbern24a15} shows, the average local wave number partially mimics the pulse number.

\begin{figure}[t]
\begin{subfigure}{.49\textwidth}
  \centering
 		\def\svgwidth{\columnwidth}
    		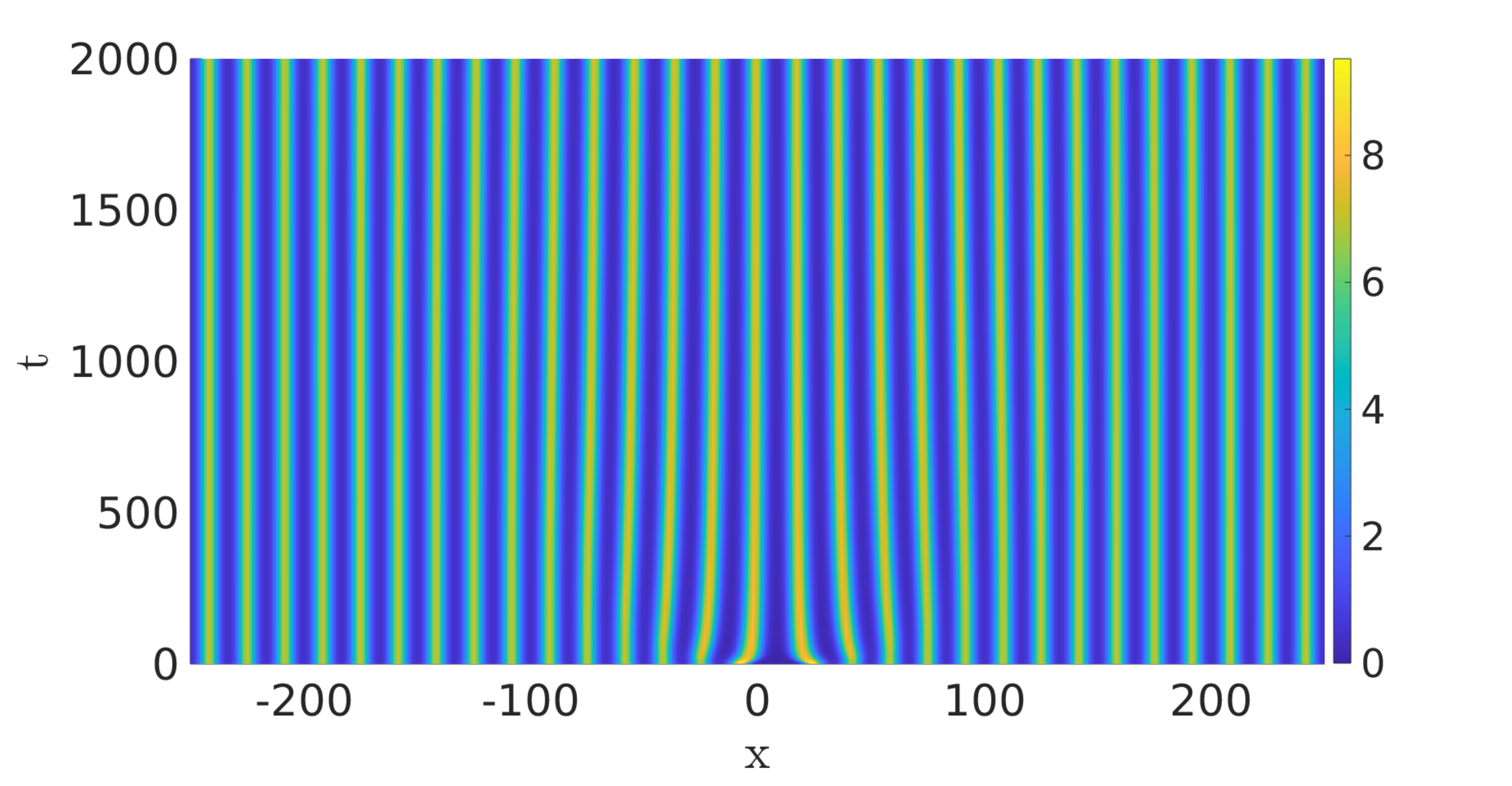
  \caption{}
    \label{fig:DetLocalwnV}
\end{subfigure}
\begin{subfigure}{.49\textwidth}
  \centering
 		\def\svgwidth{\columnwidth}
    		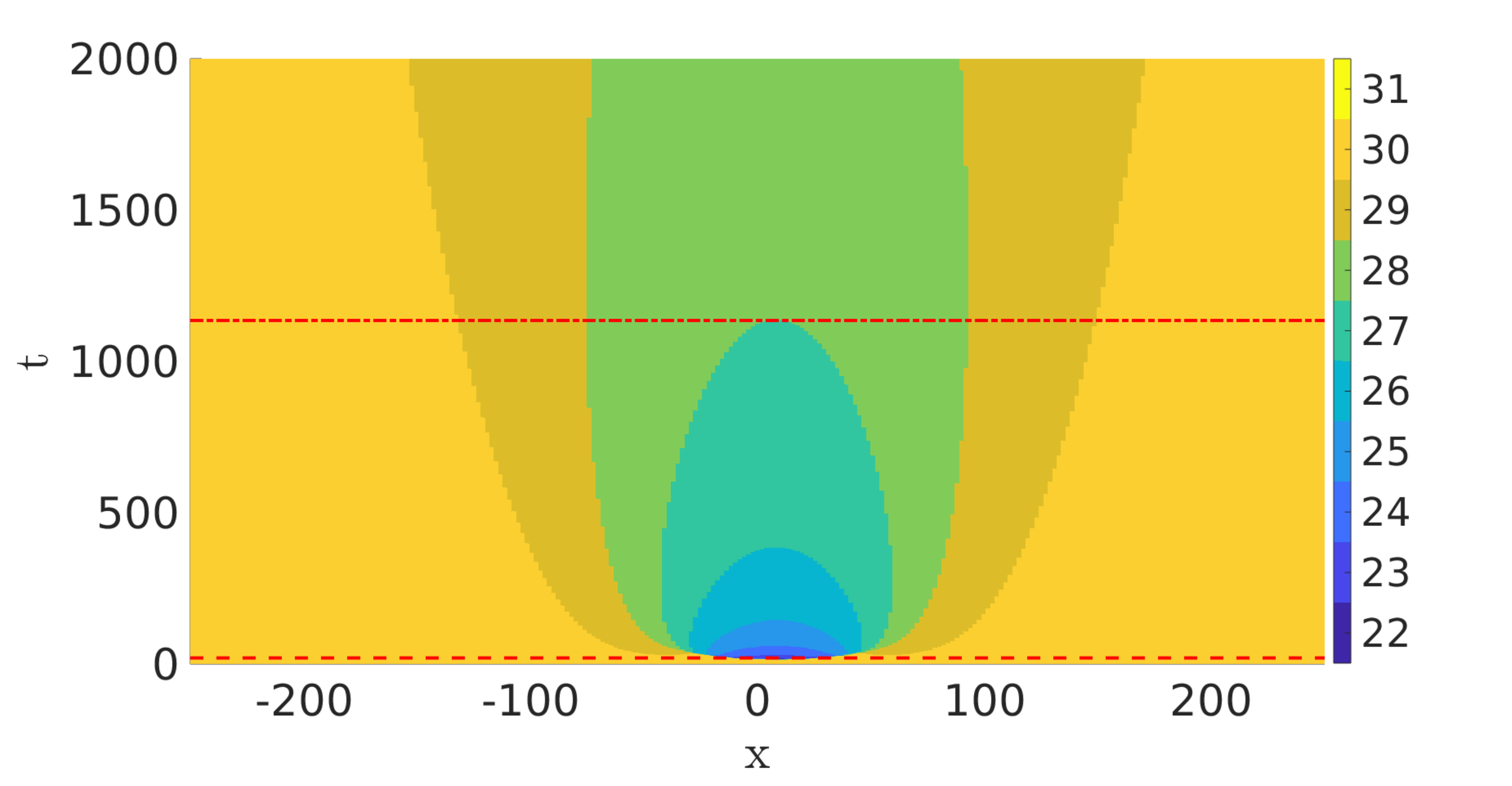
  \caption{}
    \label{fig:DetLocalwn}
\end{subfigure}
\caption{Simulation of \eqref{eq:int:StochK} with $\sigma=0$ and $a=1.5$. As initial condition we took the stationary solution with wave number $30$, but set it to $0$ on the interval $[0,\lambda]$ for both components, where $\lambda$ is the wavelength of the periodic pattern. Figure~(a) shows that the remaining $29$ pulses fill the gap created at $t=0$. In Figure~(b), the local wave number with $\ell=50$ is shown, which indicates that the spreading speed of this perturbation is slow. Even though the plot on the left seems to indicate that the solution quickly converges, the pulses are not rearranged yet to a stationary distribution on this timescale. The dashed red line at $t\approx20$ indicates when the average local wave number shifts from $30$ to $29$. The dashed-dotted line at $t\approx1135$ indicates the moment the predominant wave number shifts from $30$ to $29$.}
\label{fig:DetLocal}
\end{figure}

To further explain the notion of (typical) local wave number, we plot in Figure~\ref{fig:DetLocalwnV} a deterministic solution for an initial condition with wave number $30$, but one pulse at $x=0$ is manually deleted. As the deterministic RDE simulation (so $\sigma=0$ in \eqref{eq:int:StochK}) shows, the neighbouring pulses quickly move into the gap and the pulses rearrange into what appears to be a periodic pattern with wave number~$29$. However, the local wave number of the simulation with $\ell=50$ shown in Figure~\ref{fig:DetLocalwn} indicates that the speed with which the pulses rearrange into 29 equally spread out pulses is slow and the rearranging has not completed yet at $T=2000$. 
By construction, the pulse number is $29$ throughout the simulation. In contrast, the predominant wave number only changes to $29$ after a long ($t\approx 1135$) integration time (indicated by the dashed-dotted red line in Figure~\ref{fig:DetLocalwn}). In other words, the solution is in a transient phase until $t\approx 1135$. The average local wave number in this case sits in between the pulse and predominant wave number. It still assigns $30$ to the initial dynamics, but shifts quickly to $29$, as indicated by the red dashed line in the figure, at $t\approx 20$.

\begin{figure}
\begin{subfigure}{.49\textwidth}
  \centering
 		\def\svgwidth{\columnwidth}
    		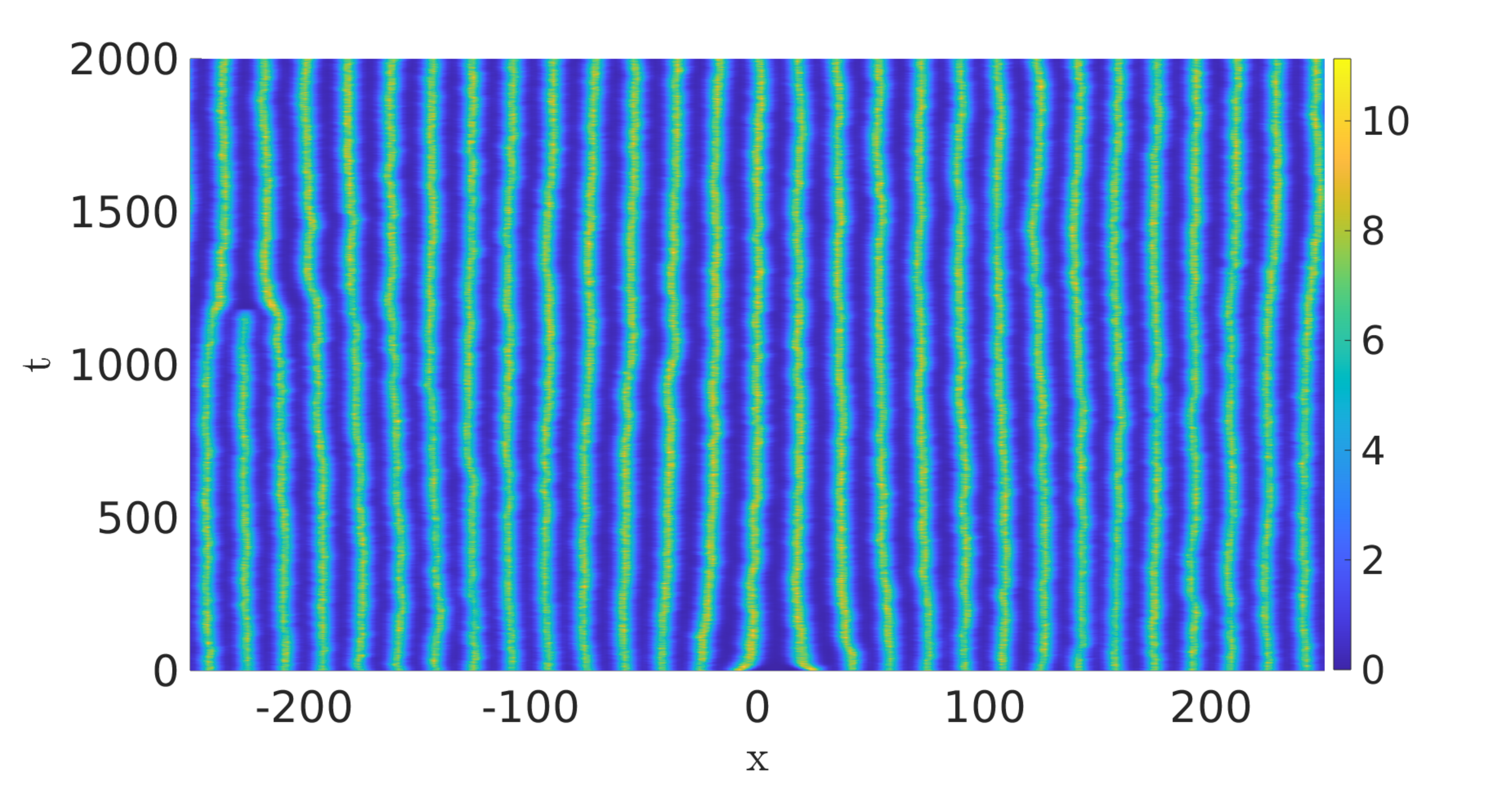
  \caption{}
    \label{fig:n30a15localwnV}
\end{subfigure}
\begin{subfigure}{.49\textwidth}
  \centering
 		\def\svgwidth{\columnwidth}
    		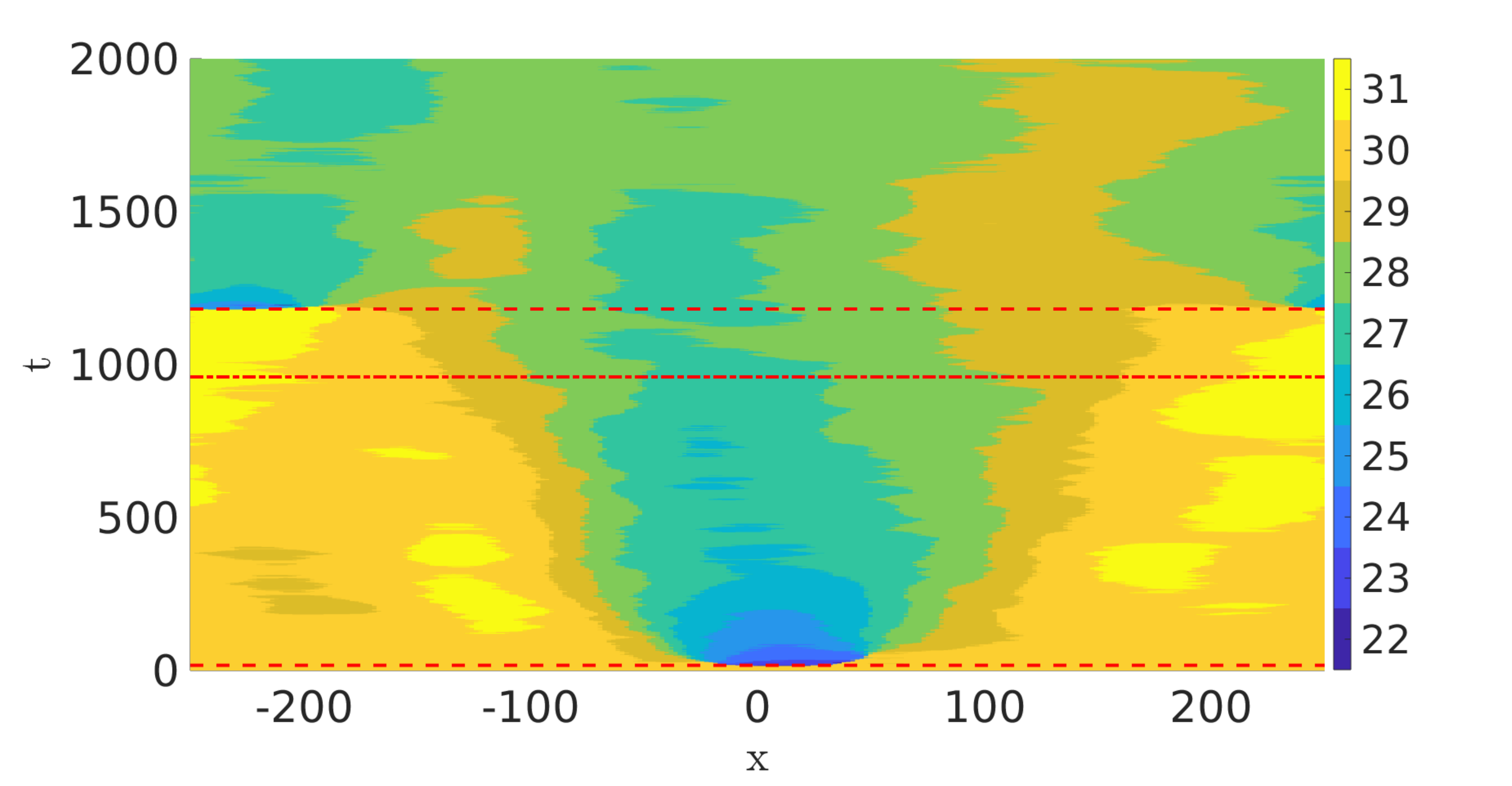
  \caption{}
    \label{fig:n30a15localwn}
\end{subfigure}
\caption{Stochastic version of Figure~\ref{fig:DetLocal} with $\sigma=0.1$. The deletion of the pulse at $x=0$ on $t=0$ is clearly visible in the local Fourier transform, as is the disappearance of a pulse around $t \approx 1200$. The two dashed lines indicate when the rounded average local wave number switches from $30$ to $29$ and from $29$ to $28$ -- the latter one coinciding with the disappearance of the pulse. Observe that just before the pulse disappears, the local wave number increases. The dashed-dotted line indicates when the predominant wave number switches from $30$ to $28$. }
\label{fig:local}
\end{figure}

In Figure~\ref{fig:local} we show the stochastic version of Figure~\ref{fig:DetLocal} with $\sigma=0.1$. As expected, it shows a perturbed version of the deterministic figure and patches with the same local wavenumber are clearly visible. The two red dashed lines indicate when the average local wave number shifts to a lower value -- and these two instances coincide with the disappearance of a pulse in the simulation. In Figure~\ref{fig:n30a15localwn}, these events are visible as a sudden decrease in the local wave number. The dashed-dotted line indicates when the predominant wave number shifts (from $30$ to $28$).

A striking difference between the deterministic and stochastic version is the presence of patches with local wave number $31$ in the stochastic case. Hence, the noise acts locally against the deterministic dynamics that pushes the pulses to an equal distribution. Note that a pulse disappears at a patch with a high local wave number. This is not just a coincidence but happens typically in these simulations. 

The simulations showcase what the local wave numbers add to the predominant wave number and/or pulse number. The pulse number tells us about the number of pulses, but not how these pulses are distributed. In contrast, the predominant wave number tells us the global distribution of the pulses, but does not tell us if all the pulses are really there. By using the local wave number, we both get a description of the distribution of the pulses, as well as how many are there.    

\begin{figure}
\begin{subfigure}{.49\textwidth}
  \centering
 		\def\svgwidth{\columnwidth}
    		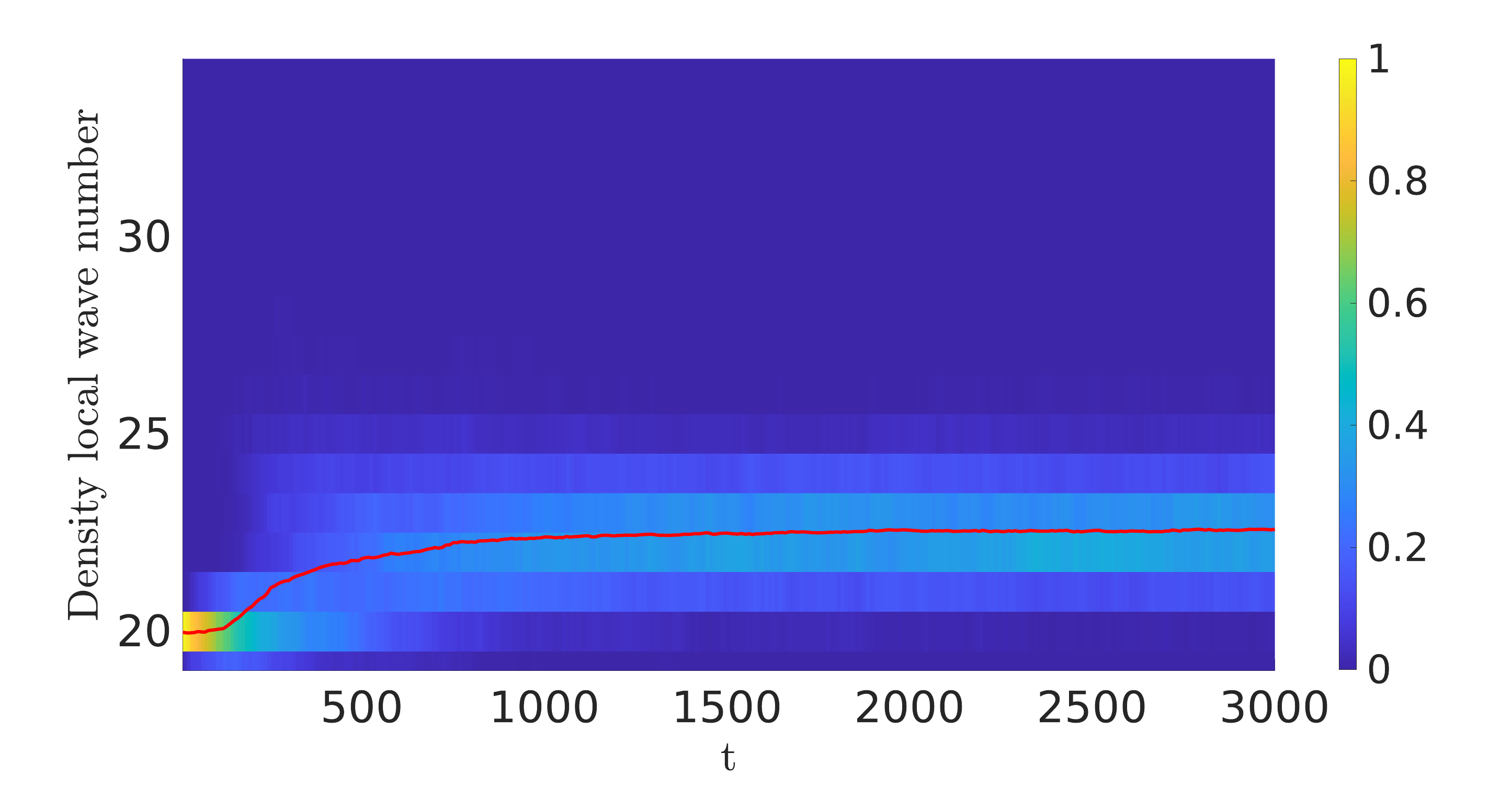
  \caption{}
    \label{fig:DensityLWn20}
\end{subfigure}
\begin{subfigure}{.49\textwidth}
  \centering
 		\def\svgwidth{\columnwidth}
    		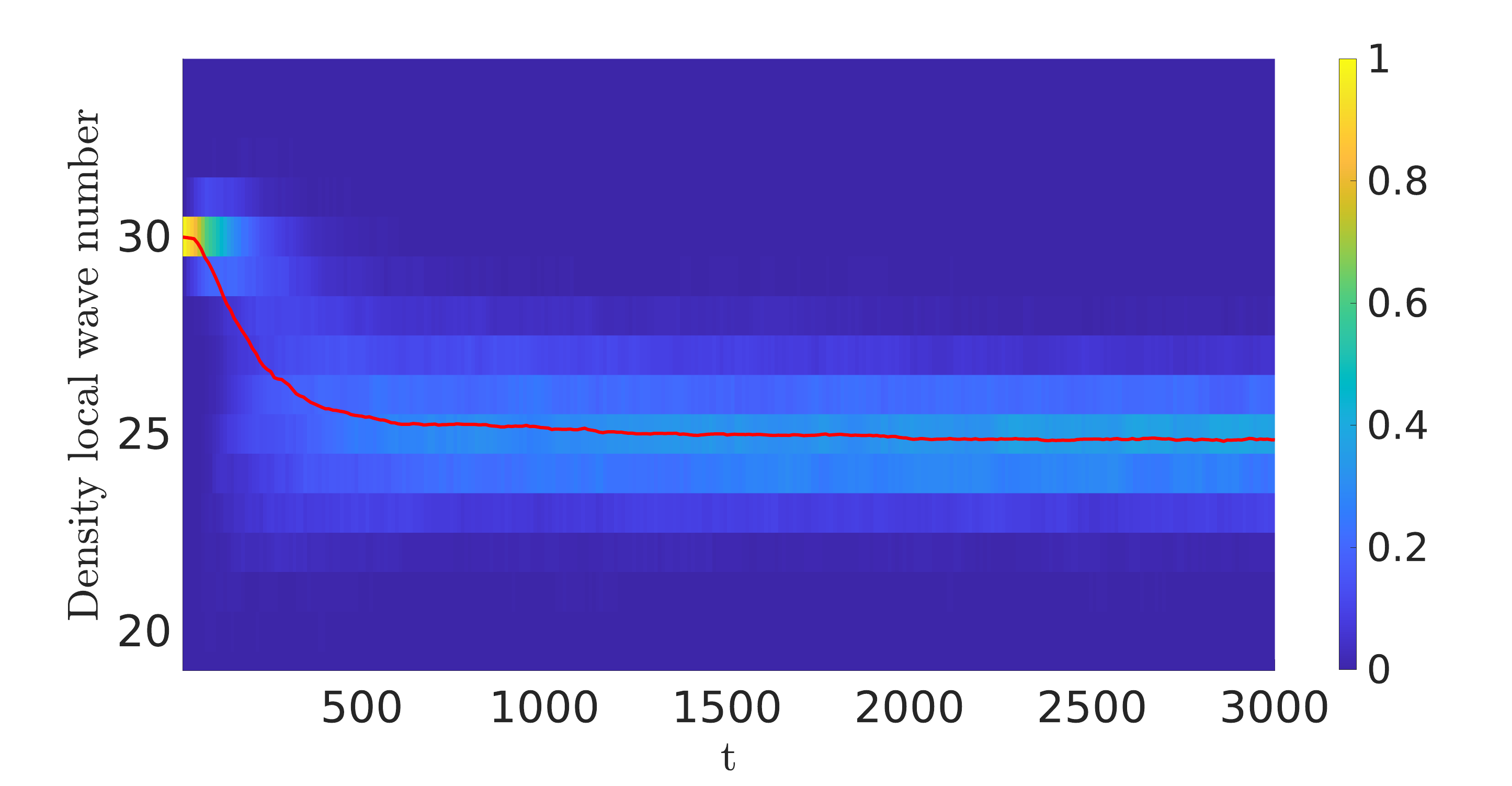
  \caption{}
    \label{fig:DensityLWn30}
\end{subfigure}
\caption{ We computed the solution to the stochastic RDE~\sref{eq:int:StochK} $50$ times with as initial condition the deterministic periodic pattern with wave number (left) $n=20$ and (right) $n=30$. We used $a=1.5$ and $\sigma=0.2$. We computed for each timestep the average distribution of local wave numbers. The red line is the average of the averaged distributions. At $t=3000$, the average distribution of local wave numbers has converged for the initial condition with wave number $n=20$, but not for $n=30$.}
\label{fig:density}
\end{figure}

With the local wave number at hand, we now have a powerful tool to track the evolution of a single realization over time. In Figure~\ref{fig:density} we show for two different initial conditions the evolution of the density of local wave numbers averaged over $50$ iterations. This means that for each timestep and iteration, we compute the local wave number for each gridpoint, which we collect into a histogram that we subsequently normalize by $2L/h$, {\emph{i.e.}} the number of gridpoints when $h$ is the spatial discretization size. These $50$ normalized histograms are then averaged resulting in an average density of local wave numbers. The average densities can be used to compute the average of the averaged density, represented by the red line in Figure~\ref{fig:density}, where it mainly serves to visualize the drift in the average density. Observe that the simulation starting from $20$ pulses quickly converges to an average of $22.5$, but when we start from $30$ pulses the average has not converged yet (over the taken simulation time) and is still decreasing. We will further study stationary distributions and these subtle differences in long-time dynamics using the averaged local wave numbers in \S\ref{sec:StatDist}.

\subsection{(Average) first exit times}
For the deterministic RDE, the (in)stability of patterns is understood intuitively. Unstable patterns will, when perturbed, transition to another (stable) pattern, while stable patterns return to the same state after a small perturbation. We generalized this idea of destabilization to the stochastic equation via the first exit time $T_\mathrm{exit}$. Upon using a stationary spatially-periodic solution with wave number $k$ as initial condition, we simulate~\eqref{eq:int:StochK} with noise~$\sigma$ and determine the first time the pattern destabilizes and call this time $T_{\rm{exit}}=T_{\rm{exit}}(a,k,\sigma,T_\mathrm{max},\omega)$, where~$\omega$ denotes the particular realization. We define this time of destabilization as the first time the solution becomes transient, see also Figure~\ref{fig:SingleR}. In practice, this implies that we measure the first time the number of pulses, {\emph{i.e.}} the pulse number, changes, as it typically changes before the predominant wave number changes. When the number of pulses does not change on the interval $[0,T_\mathrm{max}]$, we set the first exit time as~$T_\mathrm{max}$. As $T_{\rm{exit}}$ is inherently stochastic, we run many experiments to estimate the average first exit time~$E[T_{\rm{exit}}]$. Note here that determining the pulse number in a simulation is in some sense heuristic as different settings in the procedure may result in different pulse numbers.

\section{Results}
\label{results}
\subsection{(Average) first exit times}
\label{sec:ExitTime}
\begin{figure}[t]
\begin{subfigure}{.49\textwidth}
  \centering
 		\def\svgwidth{\columnwidth}
    		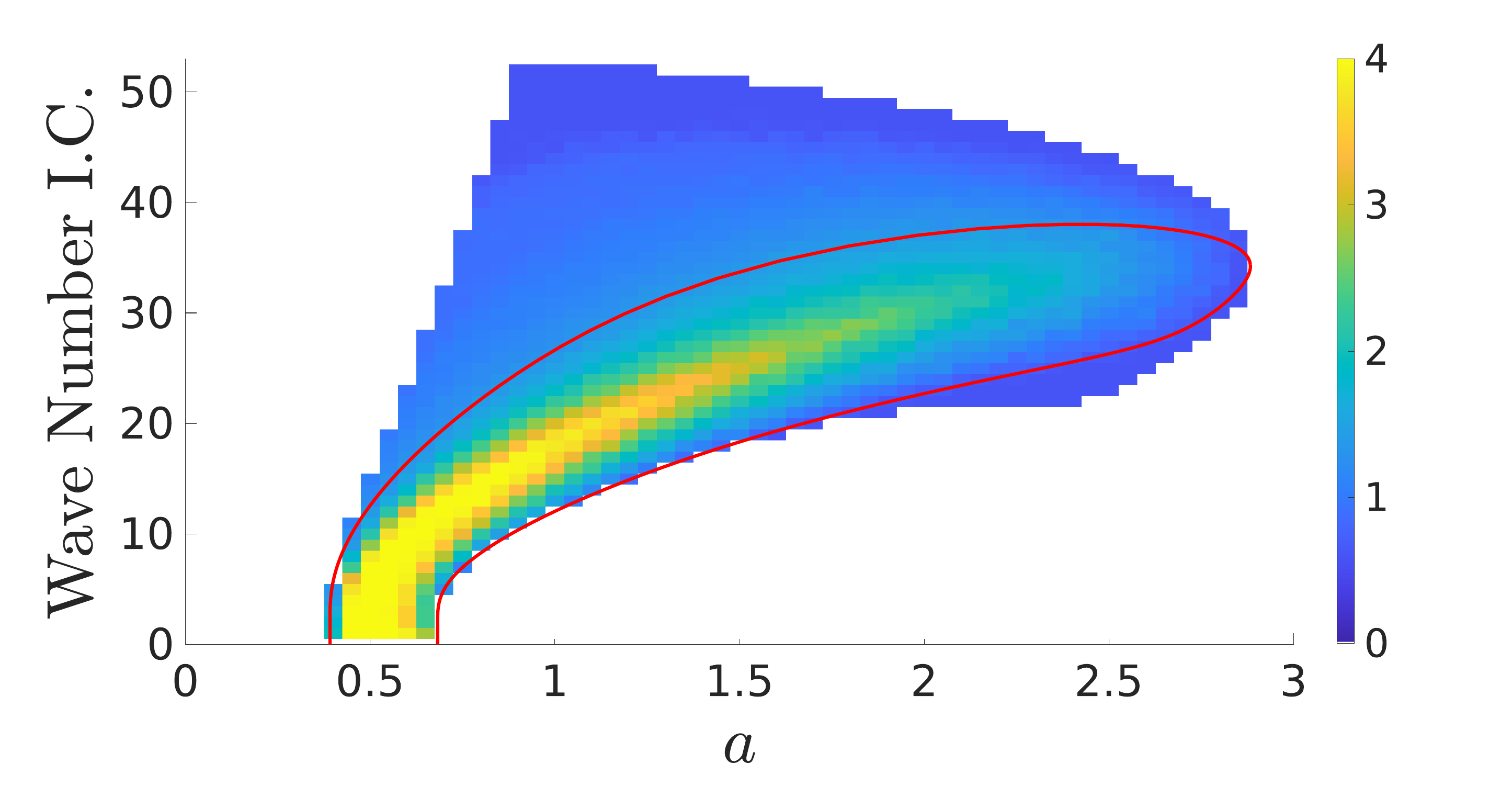
  \caption{}
    \label{fig:AFETlog}
\end{subfigure}
\begin{subfigure}{.49\textwidth}
  \centering
 		\def\svgwidth{\columnwidth}
    		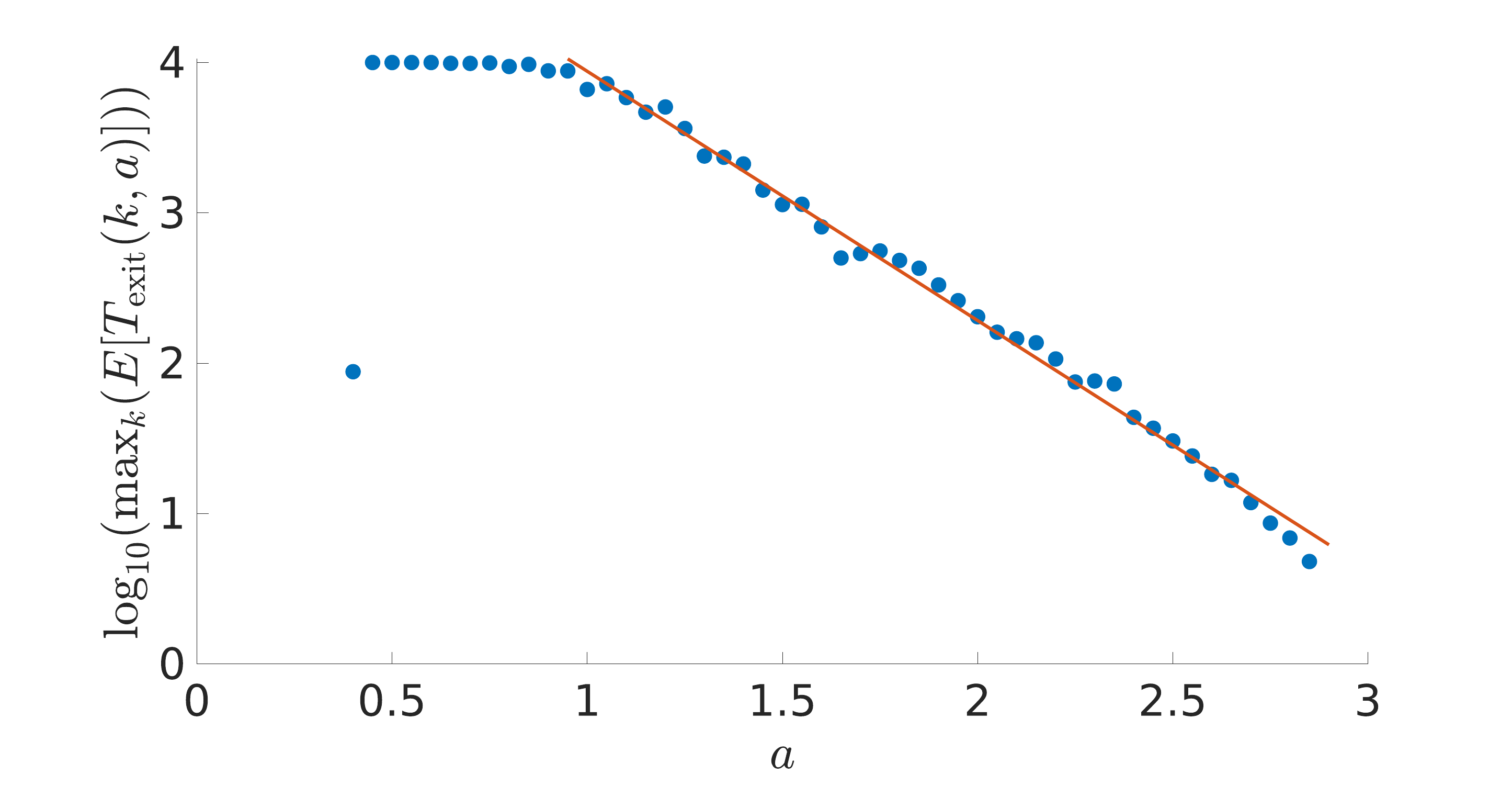
  \caption{}
    \label{fig:AFETvsA}
\end{subfigure}
\caption{Figure~(a) shows the logarithm of the average first exit time $\log_{10}(E[T_{\rm{exit}}(k,a)])$ for a range of periodic patterns used as initial conditions, and precipitation values $a$. The average is computed over 25 iterations and $\sigma=0.2$. The red line indicates the deterministic Busse balloon. Note that this figure only tells us when a pattern destabilizes on average but does not show what happens after destabilization. Figure (b) shows for each $a-$value the logarithm of the maximum average first exit time $\log_{10}(\max_kE[T_{\rm{exit}}(k,a)])$. The straight red line indicates that there is an exponential relation between $a$ and $\max_kE[T_{\rm{exit}}(k,a)]$. The deviating value for $a=0.40$ results from the fact that all periodic patterns are deterministically unstable for this value. To speed up the computation, pulse numbers are only computed for $t$ a multiple of $4$.}
\label{fig:exittimes}
\end{figure}

Here we report numerical estimates of the average first exit time $E[T_{\rm{exit}}]$ as a function of precipitation~$a$ and wave number $k$ of the initial condition, with $\sigma$ and $T_\mathrm{max}$ fixed. Note that we do not vary the system parameters $m$ and~$d$ throughout the article and hence we do not study the influence of the parameters $m$ and~$d$ on $T_{\rm{exit}}$. The result can be seen in Figure~\ref{fig:AFETlog}. As the values of $E[T_{\rm{exit}}]$ are spread over several orders of magnitude, we plot $\log_{10}\left(E[T_{\rm{exit}}]\right)$.
It is clear that for a fixed $a$-value the patterns in the middle of the Busse balloon have a significantly longer average exit time. Furthermore, observe that for some values in the centre of the balloon $E[T_{\rm{exit}}]$ is not properly resolved as many realizations did not destabilize before the end of the simulation at $T_\mathrm{max}=10^4$.

As Figure~\ref{fig:AFETlog} shows, the noise introduces a timescale for destabilization, which heavily depends on both $a$ and the position of the wave number $k$ of the initial condition in the balloon. Given the long exit times involved, it is not feasible to numerically approximate $E[T_{\rm{exit}}]$ as a function of $a, k$ and $\sigma$ simultaneously. 
From the data in Figure~\ref{fig:AFETlog}, it is however possible to directly plot~$\log_{10}(\max_kE[T_{\rm{exit}}(k,a)])$ versus $a$ with $\sigma$ fixed, as shown in Figure~\ref{fig:AFETvsA}. This shows that there is an exponential relation between $a$ and the maximum average first exit time, at least in the parameter region where we can properly resolve the average first exit time.  Similarly, as discussed in Fig.~\ref{fig:app:AFETvsSigma}, we can fix $a$ and $k$ and vary $\sigma$, resulting in a power law between $\sigma$ and $E[T_\mathrm{exit}]$.

\begin{figure}[t]
\begin{subfigure}{.33\textwidth}
  \centering
 		\def\svgwidth{\columnwidth}
    		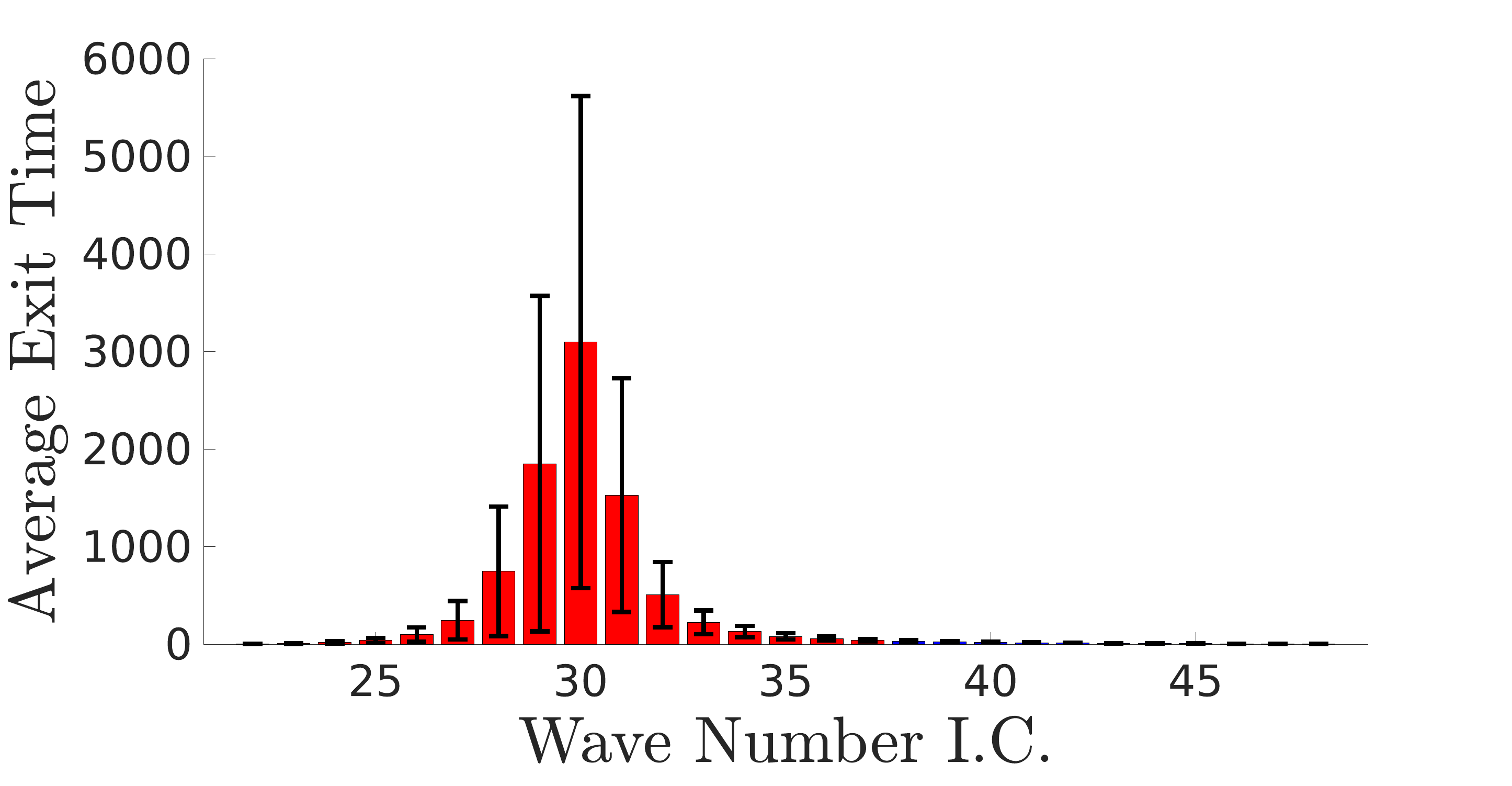
  \caption{$\sigma=0.15$}
    \label{fig:ETsigma2}
\end{subfigure}
\begin{subfigure}{.33\textwidth}
  \centering
 		\def\svgwidth{\columnwidth}
    		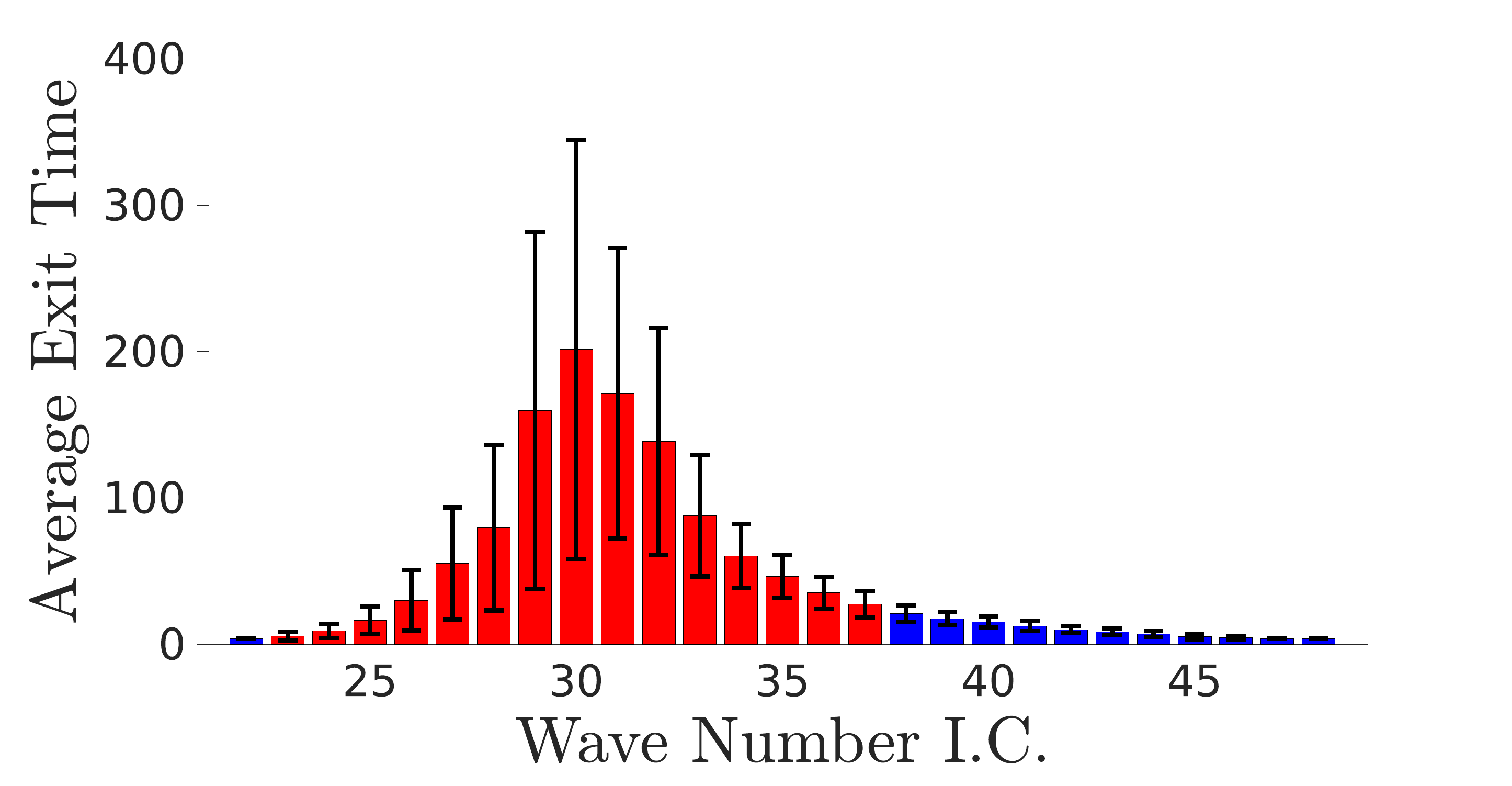
  \caption{$\sigma=0.20$}
    \label{fig:ETsigma3}
\end{subfigure}\begin{subfigure}{.33\textwidth}
  \centering
 		\def\svgwidth{\columnwidth}
    		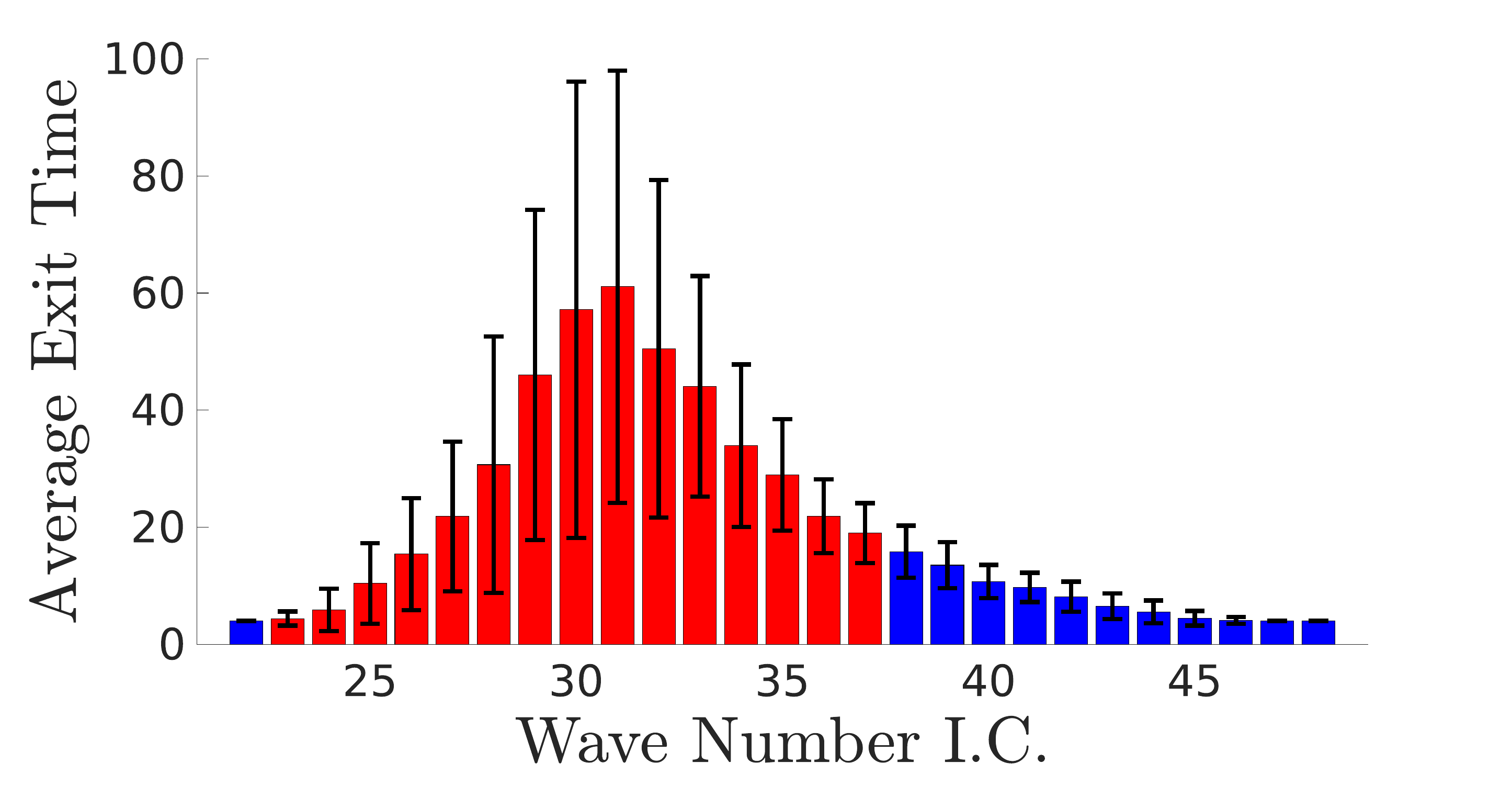
  \caption{$\sigma=0.25$}
    \label{fig:ETsigma4}
\end{subfigure}
\caption{Average exit time $E[T_\mathrm{exit}]$ for $a=2.0$ and three different values of the noise $\sigma$ computed over 50 iterations. The red bars indicate the wave numbers within the Busse balloon, and the blue ones are those outside. For each bar, the standard deviation is added to highlight the high variability of the exit times. Note that the vertical axes are different.}
\label{fig:Slices}
\end{figure}

We observe that the average first exit time depends smoothly on the wave number of the initial condition. We showcase this by plotting a single vertical slice of the Busse balloon, {\emph{i.e.}} for a fixed $a-$value and three different $\sigma$-values, see Figure~\ref{fig:Slices}. In  Figure~\ref{fig:ETsigma4}, we clearly see that the average first exit time for wave number $38$ is longer than that for $23$ and $24$, even though $23$ and $24$ are within the Busse balloon and $38$ is not. Hence there is no clear indication of the boundary of the Busse balloon. In other words, the stochasticity blurs the boundary of the deterministic Busse balloon. When the intensity of the noise is reduced, as in Figures~\ref{fig:ETsigma2} and \ref{fig:ETsigma3}, all the average exit times are increased, but significantly more for values in the middle of the balloon. Therefore, for small $\sigma-$values, the relevant timescales are set by the deterministically stable patterns. In the limit of $\sigma\to0$ and $T_{\rm max}\to\infty$, we can expect the edges of the Busse balloon to become sharp again.

\subsection{Local wave numbers and stationary distributions}
\label{sec:StatDist}
The average first exit time from the previous section gave us information on when a pattern destabilizes, but nothing on the resulting pattern after destabilization. We can classify the state after destabilization using the local wave number, but the noise will cause this classification to change continuously when using the distribution of local wave numbers, or to change discretely after a certain waiting time when using the rounded average local wave number. 

We will invest numerically if this changing of local wave numbers results in a stationary distribution of wave numbers. However, as Figure~\ref{fig:AFETlog} shows, the average exit times can be very long. Hence, simulating until the stationary distribution has been reached is not always feasible, especially for small values of $\sigma$ and $a$.
Therefore, we will focus on $a$ and $\sigma$ values for which the exit time is reasonably short. Specifically, $\sigma=0.25$ and $a\geq 1.5$. For these parameter values, we start at the unstable vegetated state $(\bar u,\bar v)$~\eqref{ss} and compute the average distribution of local wave numbers after integrating up to $T_{\rm max}= 10^4$. In Appendix~\ref{app:StatDist}, we discuss how we validated that the stationary distributions are indeed stationary. 

\begin{figure}[t]
\begin{subfigure}{.49\textwidth}
  \centering
 		\def\svgwidth{\columnwidth}
    		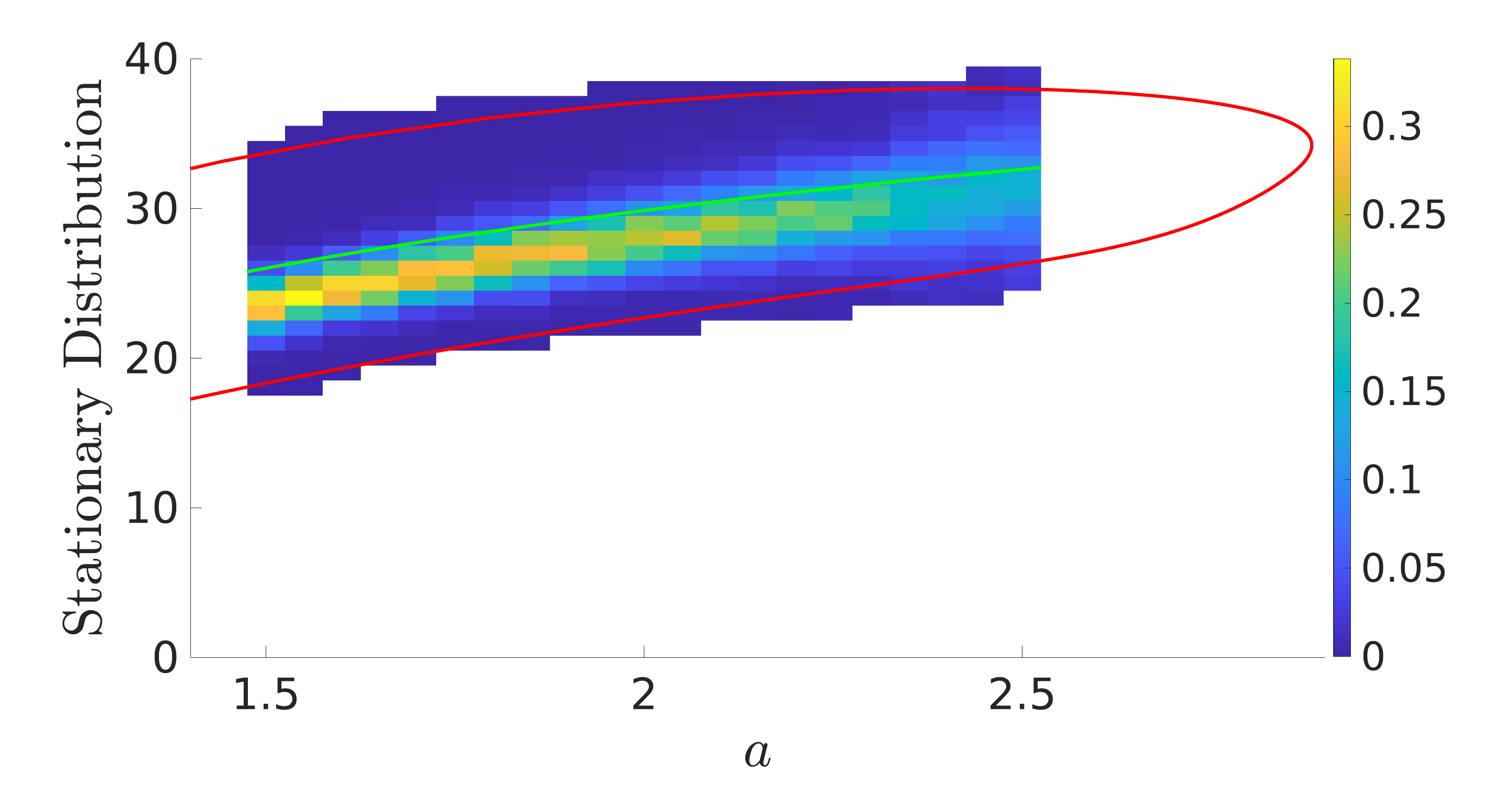
  \caption{$\sigma=0.25$}
    \label{fig:StatDistBulk}
\end{subfigure}
\begin{subfigure}{.49\textwidth}
  \centering
 		\def\svgwidth{\columnwidth}
    		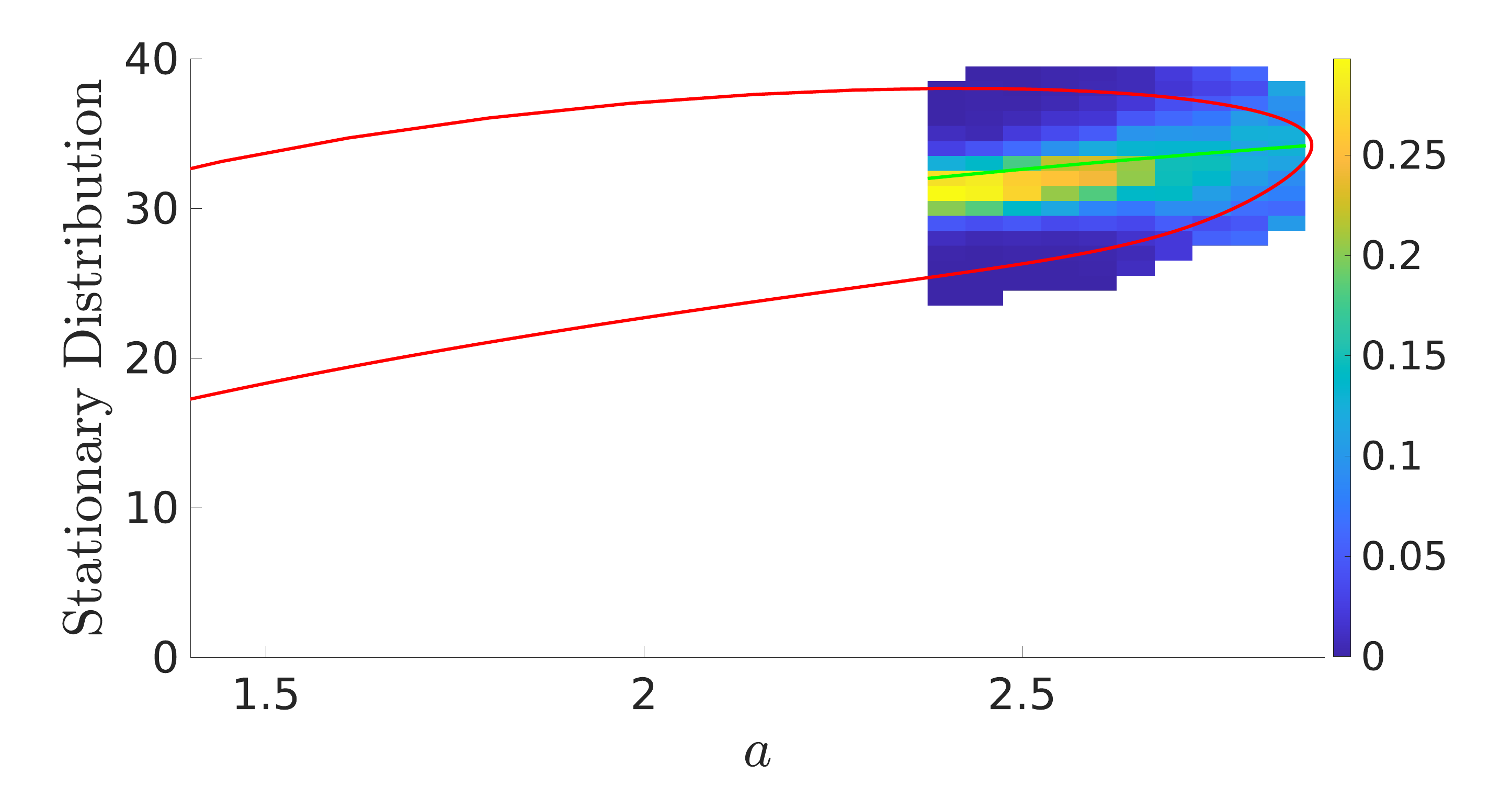
  \caption{$\sigma=0.2$}
    \label{fig:StatDistTip}
\end{subfigure}
\caption{For a section of the Busse balloon, we compute the stationary distribution of local wave numbers by starting on the unstable homogeneous state $(\bar u,\bar v)$~\eqref{ss} and integrating up to $T=T_{\rm max}=2500$ for (left) $\sigma=0.25$ and (right) $\sigma=0.2$. At the end of each simulation, the local wave numbers are counted and normalized. These normalized histograms are averaged over $100$ iterations. The red line indicates the deterministic Busse balloon and the green line the most unstable wave number~\eqref{eq:E1}.}
\label{fig:StatDist}
\end{figure}

The stationary distributions for a range of $a-$values are plotted over the Busse balloon in Figure~\ref{fig:StatDist}. In Figure~\ref{fig:StatDistBulk}, we show the stationary distributions for $\sigma=0.25$ and $a\in[1.5,2.5]$. As expected, the width of the distribution increases as $a$ increases.  For the tip of the balloon, the value of $\sigma=0.25$ is too large for meaningful patterns, see also Figure~\ref{fig:app:Sigma}. Therefore, we lower $\sigma$ to study the tip of the Busse balloon, see Figure~\ref{fig:StatDistTip}.

It must be noted that the effect the noise has on the solutions heavily depends on both $a$ and $\sigma$, hence comparing the stationary distributions for different $a-$values but fixed $\s$, might lead to skewed ideas on the shape of the distribution. For example, when we reduce the noise level to $\sigma=0.2$ in Figure~\ref{fig:StatDistTip}, we observe that the distribution becomes more peaked again. Hence, the shape of the distribution for $a=2.5$ is not essentially different from the distribution at $a=1.5$, but the effective noise values for which these shapes are attained differ.

\subsection{Comparison with previous stability measures}
\begin{figure}[t]
\begin{subfigure}{.33\textwidth}
  \centering
 		\def\svgwidth{\columnwidth}
    		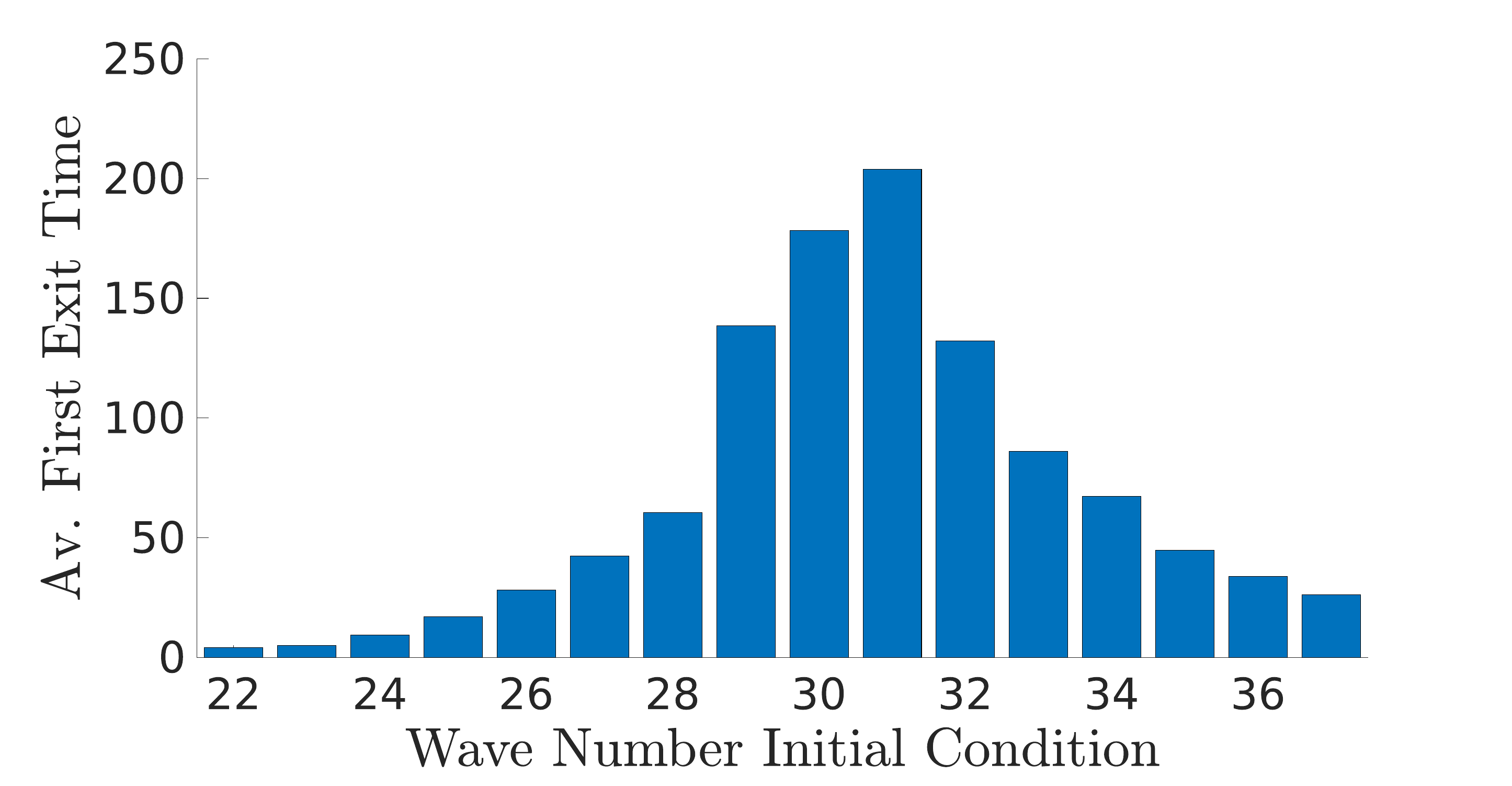
  \subcaption{Average first exit time.\; \; \; \; \; \; \; \; \; \; \; \; \; \; \; \; \; \; \; \; \; \; \; \; \; \; \; }
    \label{fig:CompAFET}
\end{subfigure}
\begin{subfigure}{.33\textwidth}
  \centering
 		\def\svgwidth{\columnwidth}
    		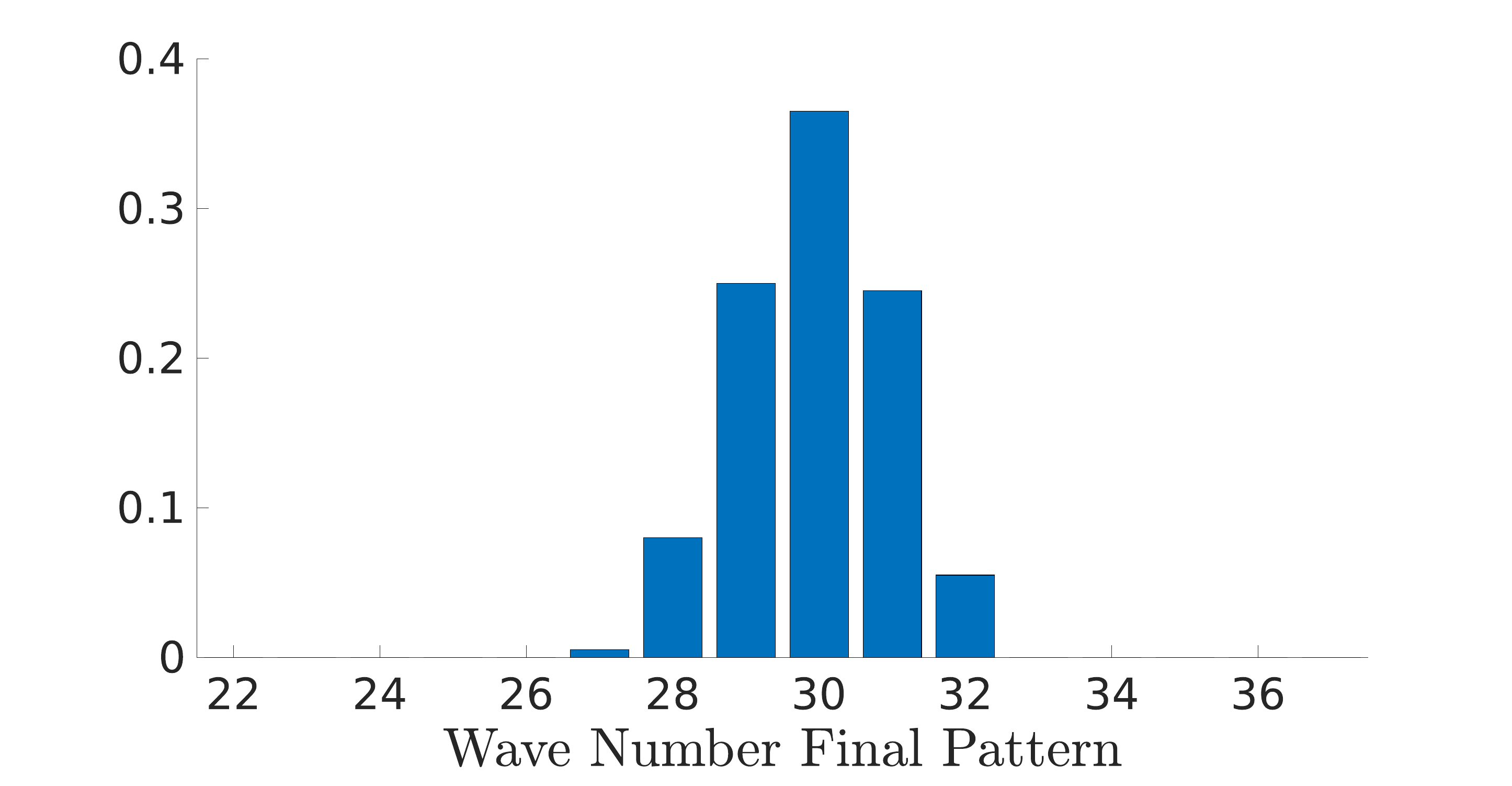
  \subcaption{Wave number from a\\\textcolor{white}{test\,}uniform background state.}
    \label{fig:CompFU}
\end{subfigure}\begin{subfigure}{.33\textwidth}
    \hspace{10mm}
  \centering
 		\def\svgwidth{\columnwidth}
    		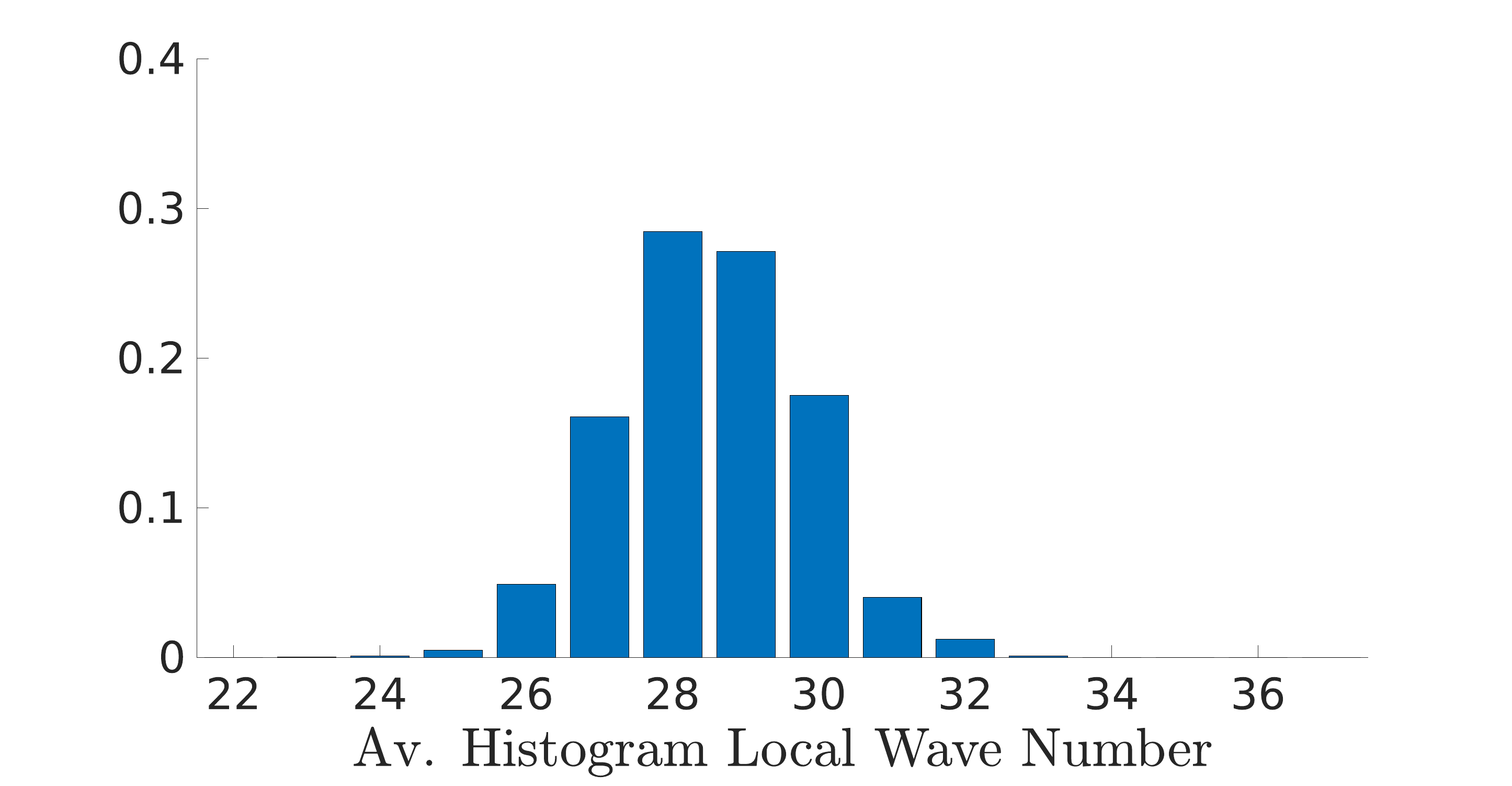
  \subcaption{Stationary Distribution. \; \; \; \; \; \; \; \; \; \; \; \; \; \; \; \; \; \; \; \; \; }
    \label{fig:CompStatDist}
\end{subfigure}
\caption{Comparison of three different approaches to stability for $a=2$ and $\sigma=0.2$. In (a) we show a slice of Figure~\ref{fig:AFETlog}, (b) shows the distribution of wave numbers in a deterministic setting such as Figure~\ref{fig:int:FromUnif}, and (c) shows a stationary distribution as in Figure~\ref{fig:StatDistBulk} but recomputed for $\sigma=0.2$.}
\label{fig:comp}
\end{figure}

With the techniques developed in the previous sections, we can now compare the different views on and measures for {\emph{(stochastic) stability}} and {\emph{observability}}. For clarity, we will focus on a single value $a=2$ and a single noise intensity $\sigma=0.2$. In Figure~\ref{fig:comp}, we show the three different approaches used: in Figure~\ref{fig:CompAFET}, we show a slice of Figure~\ref{fig:AFETlog} related to the average first exit time, Figure~\ref{fig:CompFU} shows the distribution of wave numbers in a deterministic setting such as Figure~\ref{fig:int:FromUnif}, and Figure~\ref{fig:CompStatDist} shows the stationary distribution as in Figure~\ref{fig:StatDistBulk} (but recomputed for $\sigma=0.2$). We note that the wave numbers for which the maximum is attained differ. The average first exit time attains its maximum at $31$, while this is $30$ in the deterministic experiment and $28$ for the stationary distribution. Hence, it appears that patterns typical in deterministic experiments, more often lose a pulse than gain one upon adding noise to the system.

\section{Conclusion and Discussion}
\label{sec:disc}
We have shown that for solutions of the stochastic Klausmeier model \eqref{eq:int:StochK} with fixed rainfall $a$, wave numbers are typically found in the middle of the Busse balloon, near the most unstable mode. However, the solutions are typically not well described by a single wave number, but rather by a distribution of local wave numbers. This conclusion for the stochastic model is in contrast to the deterministic model, where the initial condition determines the observed global wave number.

\begin{figure}[t]
\begin{subfigure}{.49\textwidth}
  \centering
\includegraphics[width=\columnwidth,trim={0 2.32cm 0 0},clip]{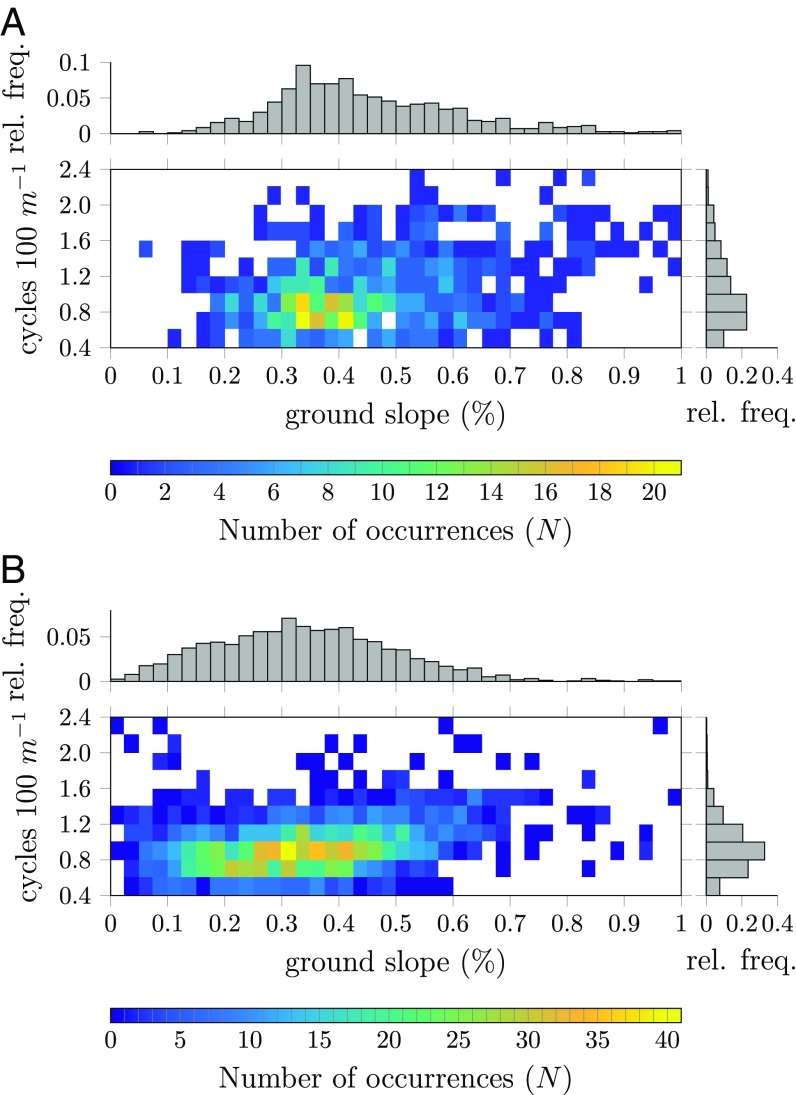}
  \caption{Site A, `Haud'}
    \label{fig:RobbinA}
\end{subfigure}
\begin{subfigure}{.49\textwidth}
  \centering
  \includegraphics[width=\columnwidth,trim={0 0 0 2.32cm},clip]{RobbinOlfa.jpg}
  \caption{Site B, `Sool'}
    \label{fig:RobbinB}
\end{subfigure}
\caption{For two sites in Somalia, Bastiaansen {\emph{et al.}}~\cite{bastiaansen2018multistability} determined the predominant wave number of stripe patterns in grids of 100 by 100 meters and determined the slope of the terrain in each grid. These figures indicate that for a single value of the slope, there is no preferred value for the predominant wave number. This can be understood as evidence for a Busse balloon. Figures reproduced form~\cite{bastiaansen2018multistability}.}
\label{fig:Robbin}
\end{figure}

The parameters $a$ and $\sigma$ in \sref{eq:int:StochK} are qualitative nondimensionalized parameters which makes it difficult to compare our results with real ecological data such as in~\cite{bastiaansen2018multistability}. However, our results could be falsified if the experimental Busse balloons would show, for example, a very uniform distribution of wave numbers or a distribution highly skewed towards the top or bottom of the balloon. As Figure~\ref{fig:Robbin} shows, this is not the case. Although the underlying model is different and the Busse balloon is studied as a function of wave number and slope instead of wave number and rainfall, these experimental Busse balloons do not contradict our results on the stochastic Busse balloon. 

To what extent can our mathematical conclusion be lifted to an ecological setting? It is tempting to conclude that also in experimental Busse balloons, we expect to observe wave numbers in the middle of the balloon. However, as we have shown, our conclusions only hold for specific combinations and ranges of $a,k,\sigma$ and~$T_\mathrm{max}$. For low $a-$values and/or low $\sigma-$values, the dynamics on long, possibly ecological-relevant, timescales, is still determined by the initial condition. Furthermore, we do not expect $a$ to remain constant on long timescales.  

In a deterministic Klausmeier model where the average rainfall $a$ slowly decreases, patterns can become unstable and have the tendency to jump back into the balloon by almost simultaneously losing half of the pulses~\cite{siteuretal2014}. In~\cite{siteuretal2014}, this scenario is replaced by adding ad hoc noise and the patterns slowly follow the upper edge of the Busse balloon, see also~\cite{Arnd2023}. Similar behaviour could be observed in \eqref{eq:int:StochK} when the rate of change $da/dt$ is faster than the typical timescale of the noise,  say $E[T_\mathrm{exit}(k,a,\sigma)]$. However, the noise will dominate the deterministic dynamics when the timescale on which the noise acts is significantly larger than this rate of change $da/dt$. This implies that we can split the dynamics into a relatively fast convergence to the stationary distribution and then a slow change in the stationary distribution as $a$ changes. Hence, for comparisons with data, it is important to be able to probe the different timescales.  

In this article, we chose explicitly to work with a reduced nondimensionalized Klausmeier model to develop the concepts of stochastic stability and observability for periodic patterns. However, there is a rich literature on more extensive deterministic models, also from different application fields~\cite[e.g]{HilleRisLambersetal2001,Giladetal2004,VdKoppeletal2005,Eppingaetal2009}, all of which could in principle be extended with noise to model inherent uncertainties in the system. Subsequently, these appended stochastic models could be numerically studied with the concepts developed in this article.

We proposed to use local wave numbers to describe the stationary distribution, but how does this compare to, say, studying the Fourier transform? In~\cite{vinals1991numerical}, the stationary distribution is studied by computing the average Fourier transform. It is reported that the width of the spectrum around the peak is wider in the stochastic case when compared with the deterministic spectrum. However, when the peak of the spectrum is located at wave number $k$, this does not imply that a typical solution can be considered as being close to the periodic pattern with wave number $k$ as we have seen. The fact that a typical solution consists of multiple patches with different wave numbers is better represented by the local wave number. 

We discussed that computing stationary distributions is not always feasible, especially for small values of $\sigma$ and $a$, because the time to convergence might be outside the range of our numerical capabilities. There are techniques that can solve this issue, {\emph{e.g.}} approaches based on large deviation theory~\cite{grafke2019numerical}, but that is not within the scope of this work. 
Another option to circumvent the problem of long integration times is to cut time into intervals of length $T$ and view the dynamics as a Markov chain on the space of wave or pulse numbers. The transition matrix $P_T$ can then be estimated just by computing the dynamics on $[0,T]$ and the dynamics for $t\to\infty$ can then be found by computing the stationary distribution of $P_T$. However, this results in two challenges. First, for numbers in the middle of the Busse balloon, the values on the diagonal of $P_T$ can be very close to $1$. In order to numerically resolve these values from $1$, either $T$ must be large or the number of iterations used to estimate $P_T$, again resulting in problems with computation time. Secondly, a more fundamental problem is that the states of the Markov chain are not well defined. Do we use pulse or wave numbers? And, assuming we choose wave numbers, what do we take as initial condition for the simulation with state $k$, do we have to generate random initial conditions from the space of all quasi-periodic patterns that have wave number $k$? These complications make the option of using Markov chains unattractive. 

For practical reasons, we restricted ourselves to one space dimension in this article. These 1D patterns correspond to a single slice of 2D stripes. This implies that the noise does not vary along the stripes, and some care should be taken when interpreting the current results in 2D as spotted or hexagonal patterns cannot be studied in one space dimension. Although the Busse balloon has some technical challenges in 2D that have not been extensively studied, the numerical techniques outlined here can be readily extended to two spatial dimensions.

\section*{Acknowledgements}
We like to thank the Wageningen University \& Research for supporting this work through the research programme Data Driven Discoveries in a Changing Climate (D3C2). The work of CH is funded by the Dutch Institute for Emergent Phenomena (DIEP) at the University of Amsterdam via the program Foundations and Applications of Emergence (FAEME) 

\section{Data availability}
All data for the figures can be generated using the scripts on \url{https://github.com/chshamster/StochasticBusseBalloon}
\bibliographystyle{klunumHJ}
\bibliography{ref}

\appendix
\renewcommand\thefigure{\thesection.\arabic{figure}}    
\setcounter{figure}{0}    
\section{Appendix}

\subsection{Computation of the most unstable mode}
\label{app:MostUnst}
Here we compute the most unstable mode of the PDE dynamics around the unstable homogeneous state $(\bar u,\bar v)=(a/2-\sqrt{a^2-4m^2}/2,(a-\bar u)/m)$. Upon linearizing the deterministic version of \sref{eq:int:StochK} around the homogeneous steady state $(\bar u,\bar v)$, the linear operator $\mathcal{L}$ becomes 
\begin{align*}
    \mathcal{L}=\begin{pmatrix}d_1\p_{xx}&0\\0&d_2\p_{xx} \end{pmatrix}+\begin{pmatrix}
        -1-{\bar v}^2 & -2\bar u\bar v \\ {\bar v}^2 & -m+2\bar u\bar v
    \end{pmatrix}.
\end{align*}
In the main text, $d_1=d$ and $d_2=1$, but we keep both here for convenience. We denote the second matrix by $A$. Then, the dispersion relation $\lambda(k)$~\cite{vdStelt2013rise} is given by
\begin{align}
    \det\left[\begin{pmatrix}-d_1k^2&0\\0&-d_2k^2 \end{pmatrix}+A-\lambda\right]=0,
\end{align}
which results in 
\begin{equation} 
\label{eq:E1}
    \lambda^2+\left((d_1+d_2)k^2-\tr(A)\right)\lambda+d_1d_2k^4-\Gamma k^2+\det(A)=0
\end{equation} where 
$\Gamma=d_1A_{22}+d_2A_{11}$. Looking for extrema of \eqref{eq:E1}, implicit differentiation with respect to $k$ yields
\begin{equation*}
    2\lambda\left(\frac{\partial\lambda}{\partial k}+(d_1+d_2)k\right)+\left((d_1+d_2)k^2-\tr(A)\right)\frac{\partial\lambda}{\partial k}+4d_1d_2k^3-2\Gamma k=0.
\end{equation*} Imposing criticality, $\frac{\partial\lambda}{\partial k}=0$, and $\lambda$ real, we find that at the most unstable mode $k_\mathrm{mu}$, the corresponding eigenvalue is given by 
\begin{equation} \label{eq:E2}
    \lambda_{\mathrm{mu}}=\frac{\Gamma-2d_1d_2k^2}{d_1+d_2}.
\end{equation}
Substitution of \eqref{eq:E2} into \eqref{eq:E1} results in an even fourth-order polynomial in $k_\mathrm{mu}$. For the values of $d_1,d_2$ and $m$ chosen in the main text, the branch of $k_\mathrm{mu}=k_\mathrm{mu}(a)$ in the positive quadrant coincides with the most unstable mode. However, note that the value of $\lambda$ attained at $k_\mathrm{mu}$ is not necessarily positive. In fact, it changes sign exactly at $a=a_T$, {\emph{i.e.}} at the Turing bifurcation. In the figures in the main text, the branch $k_\mathrm{mu}$ is drawn using the implicit plot functions in Matlab.

\subsection{Numerical integration scheme}
\label{app:NumScheme}
The numerical integration scheme that we use to simulate \sref{eq:int:StochK} is based on~\cite{lord2014book}. 
When we write $X(t)$ for $(u(t),v(t))^T$, $L$ for the linear part of \sref{eq:int:StochK}, $f$ for the nonlinear part and $\sigma$ for the multiplicative noise term, we get
\begin{align*}
    dX=[LX+f(X)]dt+\sigma (0,v)^TdW^Q_t.
\end{align*}
We discretize this equation as
\begin{align*}
    X(t+dt)=X(t)+dtL_hX(t+dt)+dtf(X(t))+\sqrt{dt}\sigma (0,v(t))^TdW^Q_t,
\end{align*}
where $L_h$ stands for the spatially discretized version of the linear operator, {\emph{i.e.}} the second derivate is approximated using the standard second-order midpoint rule. Now we compute $X(t+dt)$ by solving
\begin{align}
    (I-dtL)X(t+dt)=X(t)+dtf(X(t))+\sqrt{dt}\sigma (0,v(t))^TdW^Q_t.
\end{align}
The stochastic term $dW_t^Q$ is generated using Algorithm 6.7 from~\cite{lord2014book}. The covariance function $q$ is chosen as $\frac{1}{2\xi}e^{-\pi x^2/(4\xi^2)}$, where the correlation length is chosen as $\xi=0.1$. 

\subsection{Computing wave numbers and number of pulses}
\label{app:numbers}
Numerically computing the predominant wave number in space is straightforward using the Fourier series. Computing the number of pulses is significantly more involved and in essence heuristic. First, we smooth the solution $u$ using the Matlab command \texttt{smoothdata(u,`gaussian',64)}. Note that the correlation length $64$ means $64$ gridpoints and hence is intrinsically linked to the spatial discretization $h$. Typically, we choose 64 to correspond with $h=2L/2^{12}\approx0.1221$. Once the solution is smoothened, we compute the number of local extrema using \texttt{islocalmax} and \texttt{islocalmin} with the option \texttt{MinProminence} set to 0.3. When the number of local minima and maxima is unequal, we choose the maximum of the two. On a periodic domain, we expect the two to be equal, but \texttt{islocalmax} cannot deal with periodic boundary conditions.  

The choice of smoothing parameters inherently creates a cut-off between patches that we count and those we do not. This also implies that a pattern can shortly drop above or below the threshold, resulting in a noisy signal in time while such a variation in the number of pulses is not visible in the full stochastic RDE simulation. To filter these short-time fluctuations, a moving median in time is applied.

\subsection{Stationary distribution of local wave numbers}
\label{app:StatDist}
\begin{figure}[t]
\begin{subfigure}{.49\textwidth}
  \centering
 		\def\svgwidth{\columnwidth}
    		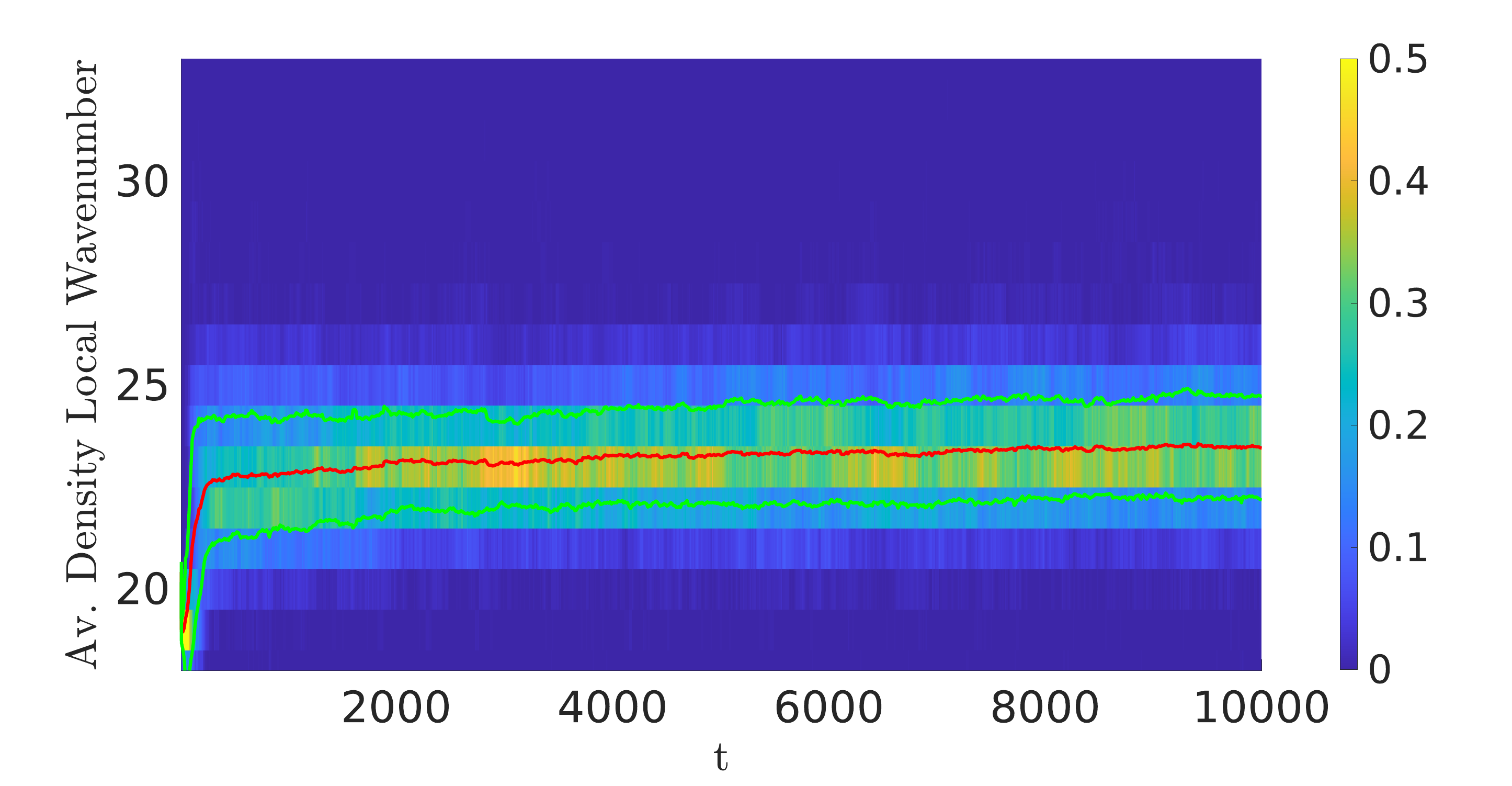
  \caption{}
    \label{fig:StatDista_5sigma25}
\end{subfigure}
\begin{subfigure}{.49\textwidth}
  \centering
 		\def\svgwidth{\columnwidth}
    		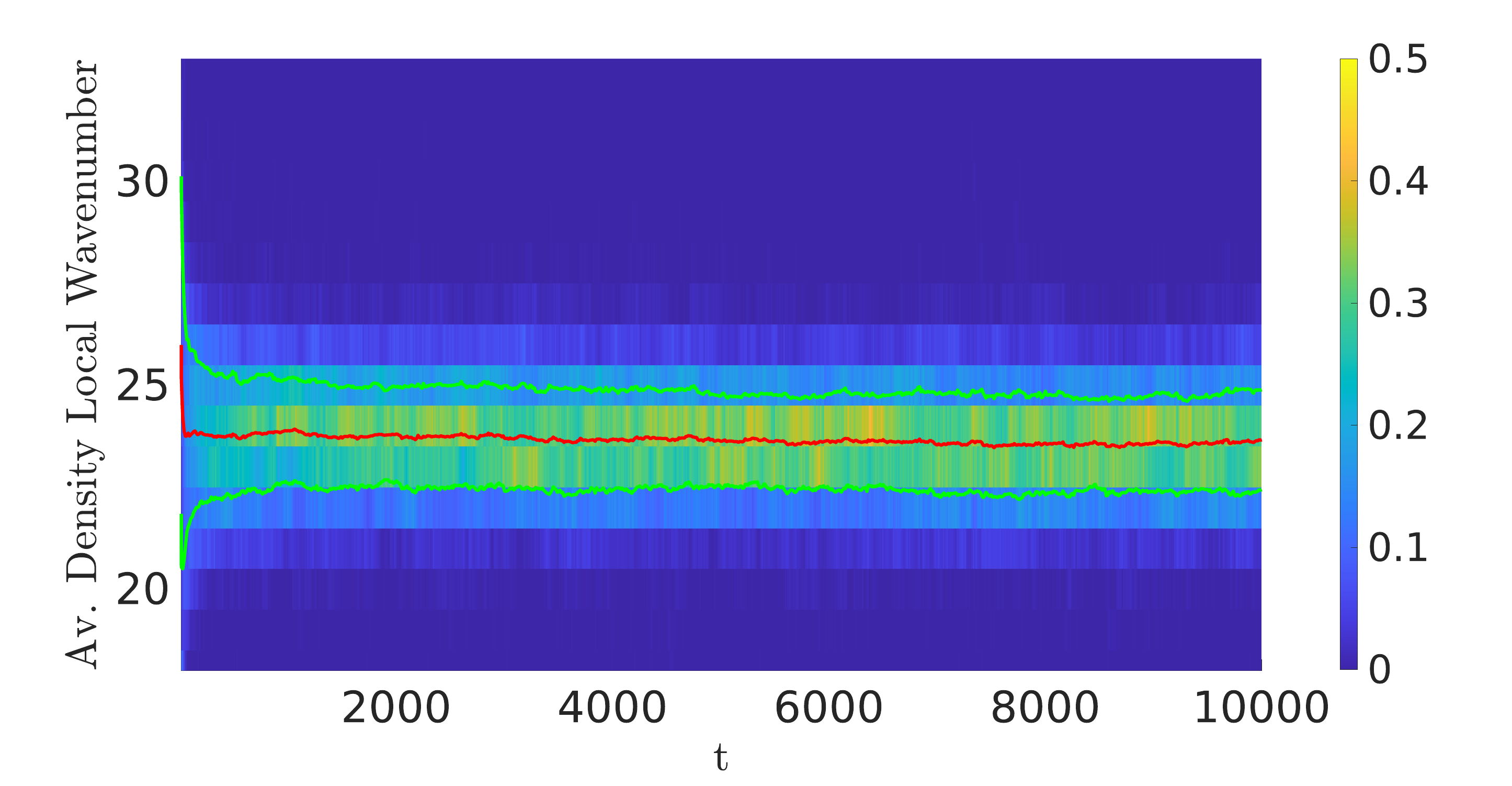
  \caption{}
    \label{fig:StatDista_5sigma25FU}
\end{subfigure}
\caption{Average density of local wave number with two different initial conditions for $a=1.5$ and $\sigma=0.25$. The average is taken over $25$ iterations. In~(a), the initial condition is a stable periodic pattern with wave number $19$, while in~(b) the initial condition is given by the unstable vegetated state $(\bar u,\bar v)$~\eqref{ss}. The red line indicates the average of the averaged densities, and the green lines indicate the standard deviation.}
\label{fig:StatDist25}
\end{figure}

One has to be careful with numerically computing stationary distributions; how does one ensure that the system is indeed in a stationary state and not slowly evolving? To test that the integration time is long enough, we run the following test for $a=1.5$, {\emph{i.e.}} the slowest value in our range for~$a$, see Fig.~\ref{fig:exittimes}. We start with two different initial conditions, the uniform vegetated state $(\bar u,\bar v)$ and the periodic state with wave number $19$. As Figure~\ref{fig:StatDist25} shows, the average distributions of local wave numbers are initially very different but converge to the same average and variance.

\begin{figure}[t]
\begin{subfigure}{.49\textwidth}
  \centering
 		\def\svgwidth{\columnwidth}
    		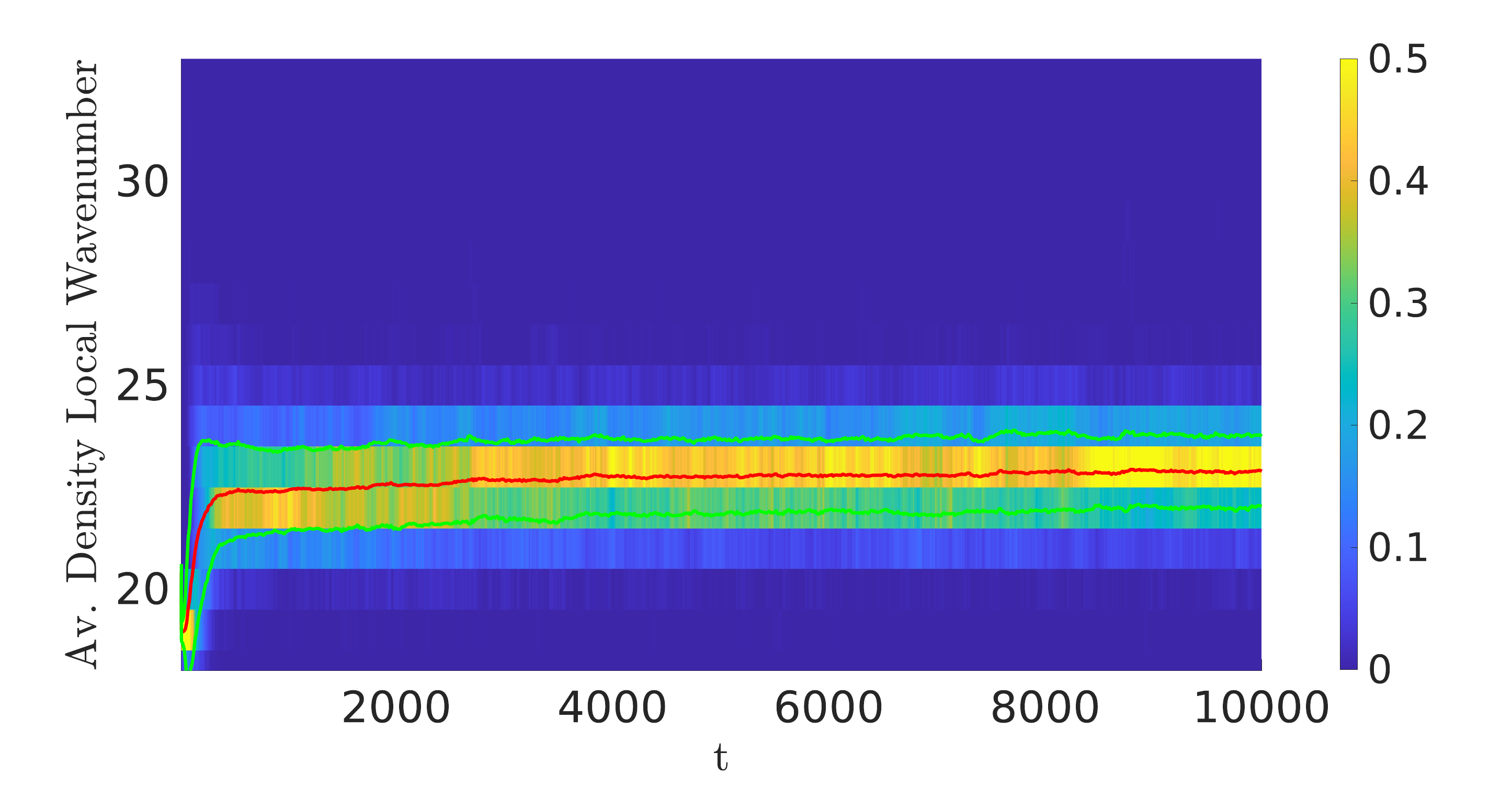
  \caption{}
    \label{fig:StatDista_5sigma20}
\end{subfigure}
\begin{subfigure}{.49\textwidth}
  \centering
 		\def\svgwidth{\columnwidth}
    		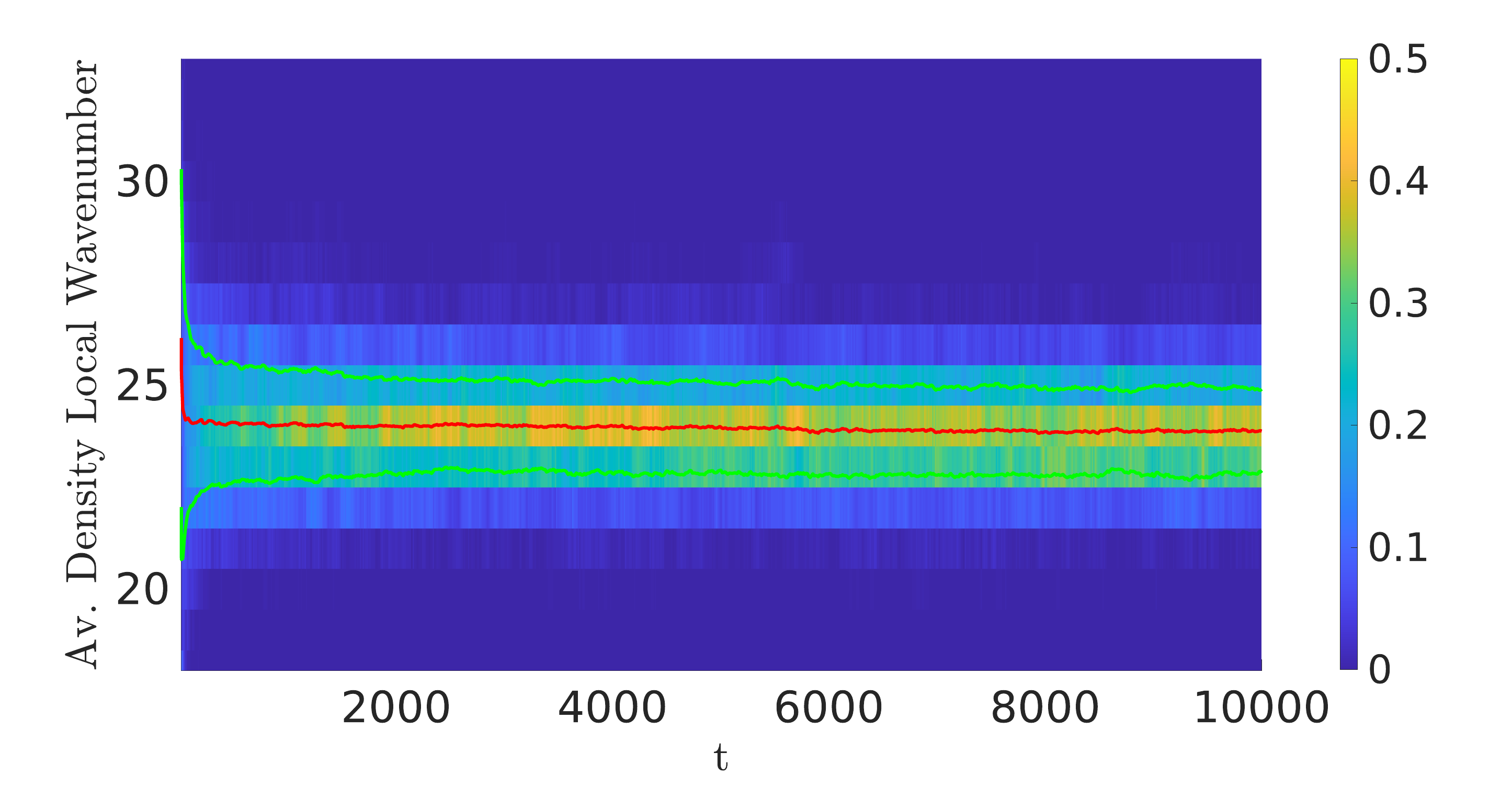
  \caption{}
    \label{fig:StatDista_5sigma20FU}
\end{subfigure}
\caption{The same as Figure~\ref{fig:StatDist25}, but with smaller noise intensity $\sigma=0.2$. Note that the distributions at $T=10^4$ are different left and right, but visually it is hard to detect that the distribution on the left is still increasing.}
\label{fig:StatDist20}
\end{figure}

When we do the same experiment but now for $\sigma=0.2$ instead of $\sigma=0.25$, we observe that the two different initial conditions lead to different distributions at $T=T_{\rm max}=10^4$, see Figure~\ref{fig:StatDist20}. That is, the simulations have not yet reached the stationary distribution. However, it is also important to note that, at least visually, the simulations seem to have converged, as both the averages and the standard deviations of the local wave number distributions did not change detectably in the final part of the simulations. As Figures~\ref{fig:StatDista_5sigma25FU} and \ref{fig:StatDista_5sigma20FU} show, the convergence to the stationary distribution is fast when the dynamics starts at $(\bar u,\bar v)$. In principle, this would allow us to use values of $a$ smaller than $1.5$, but then we would not be able to compare it with other initial conditions in order to check that we did not introduce any bias.

\subsection{Additional Figures}
\label{App:AF}
\begin{figure}[h!]
\begin{subfigure}{.49\textwidth}
  \centering
 		\def\svgwidth{\columnwidth}
    		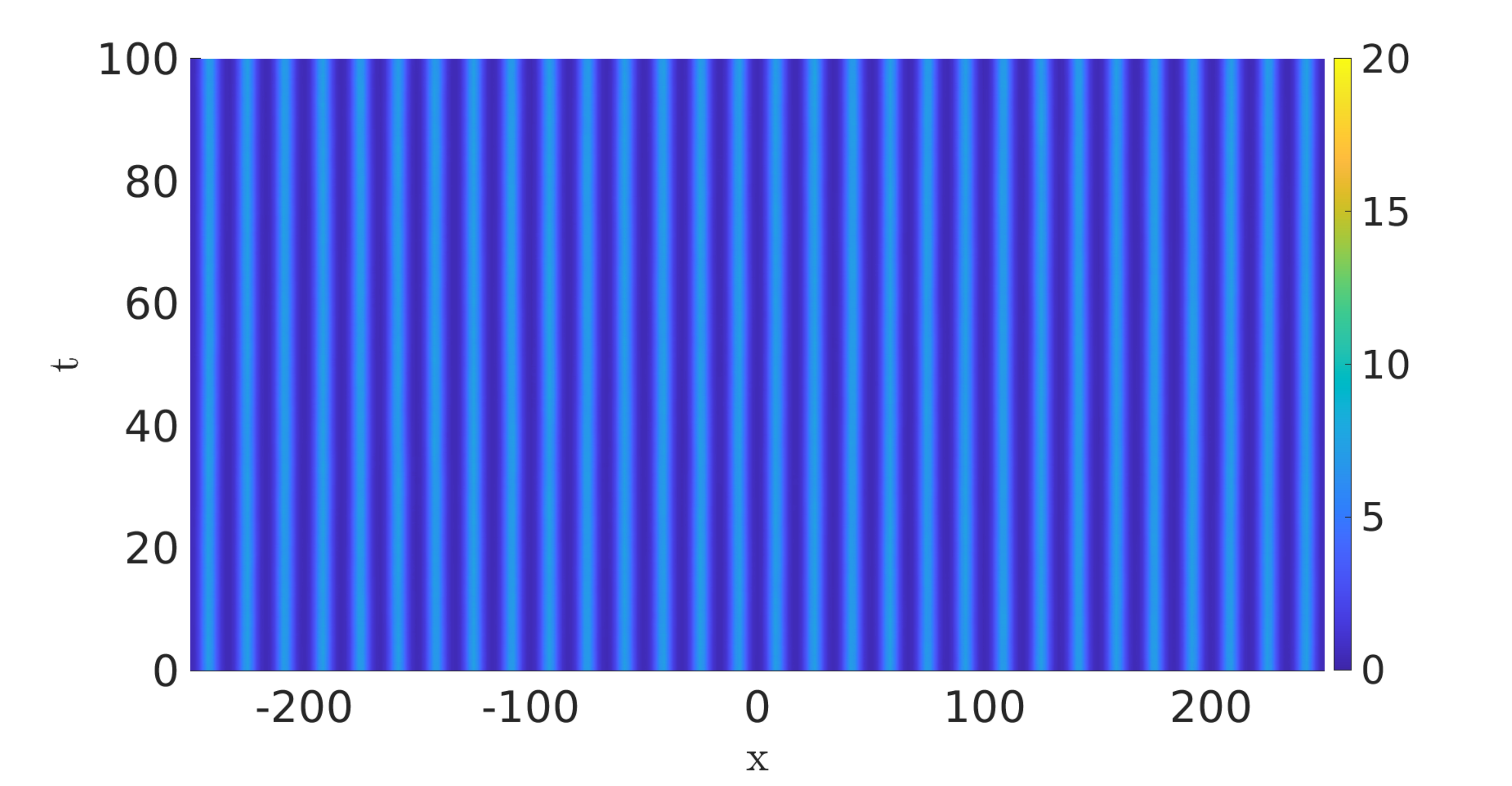
  \caption{$\sigma=0.02$}
\end{subfigure}
\begin{subfigure}{.49\textwidth}
  \centering
 		\def\svgwidth{\columnwidth}
    		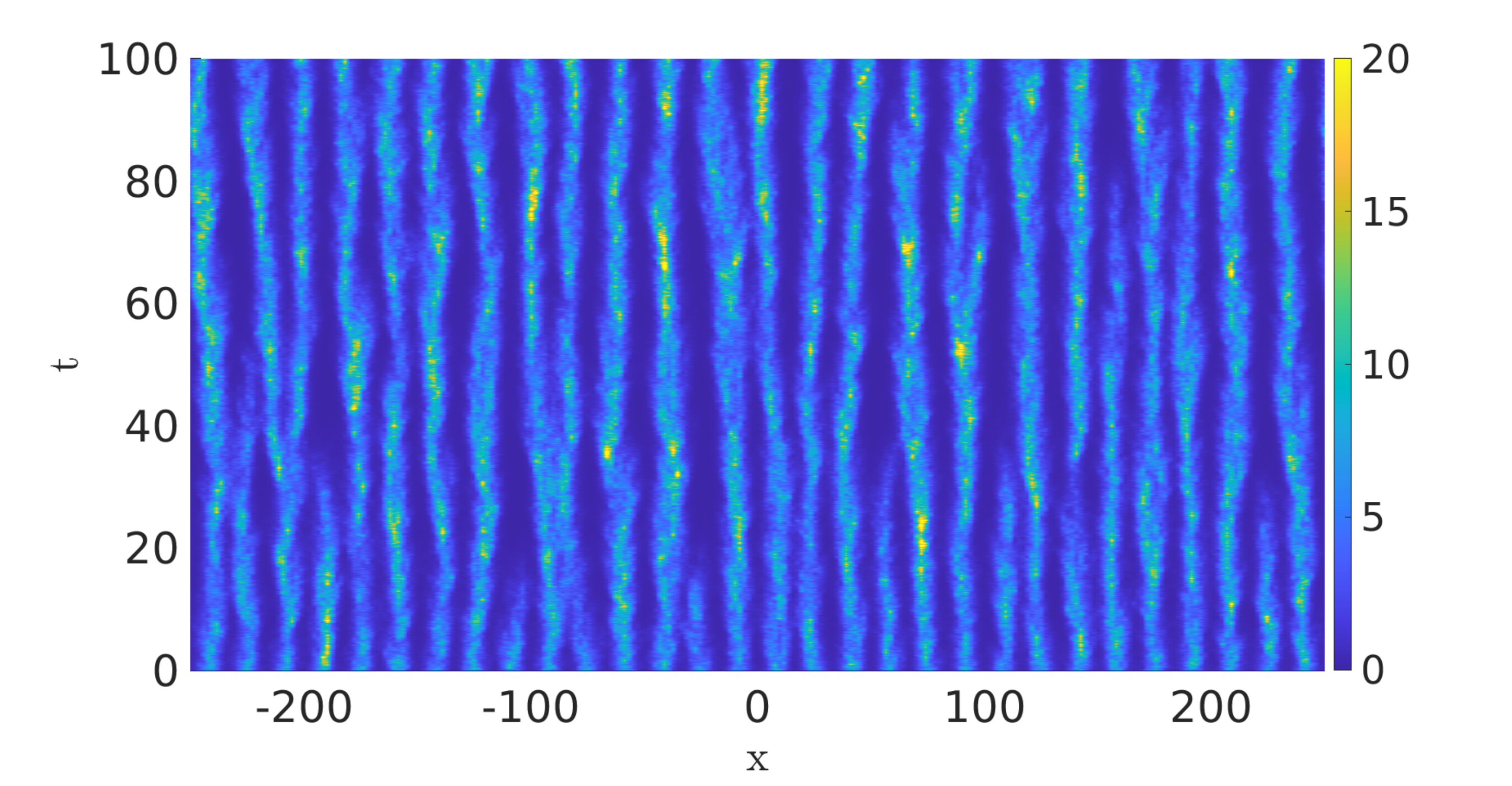
  \caption{$\sigma=0.5$}
\end{subfigure}
\caption{Two realizations of \eqref{eq:int:StochK} for $a=1.5$ and an initial condition with wave number $30$, but with different noise intensities $\sigma$. In Figure~(a), where the noise is smaller than typically used in the main text, the solution does visually not differ from the deterministic solution. In Figure~(b), the noise is larger than used in the main text and the pattern is also very irregular. Notice the bright yellow spots where the solution is significantly higher than in other places.}
\label{fig:app:Sigma}
\end{figure}

\begin{figure}
  \centering
 		\def\svgwidth{0.7\columnwidth}
    		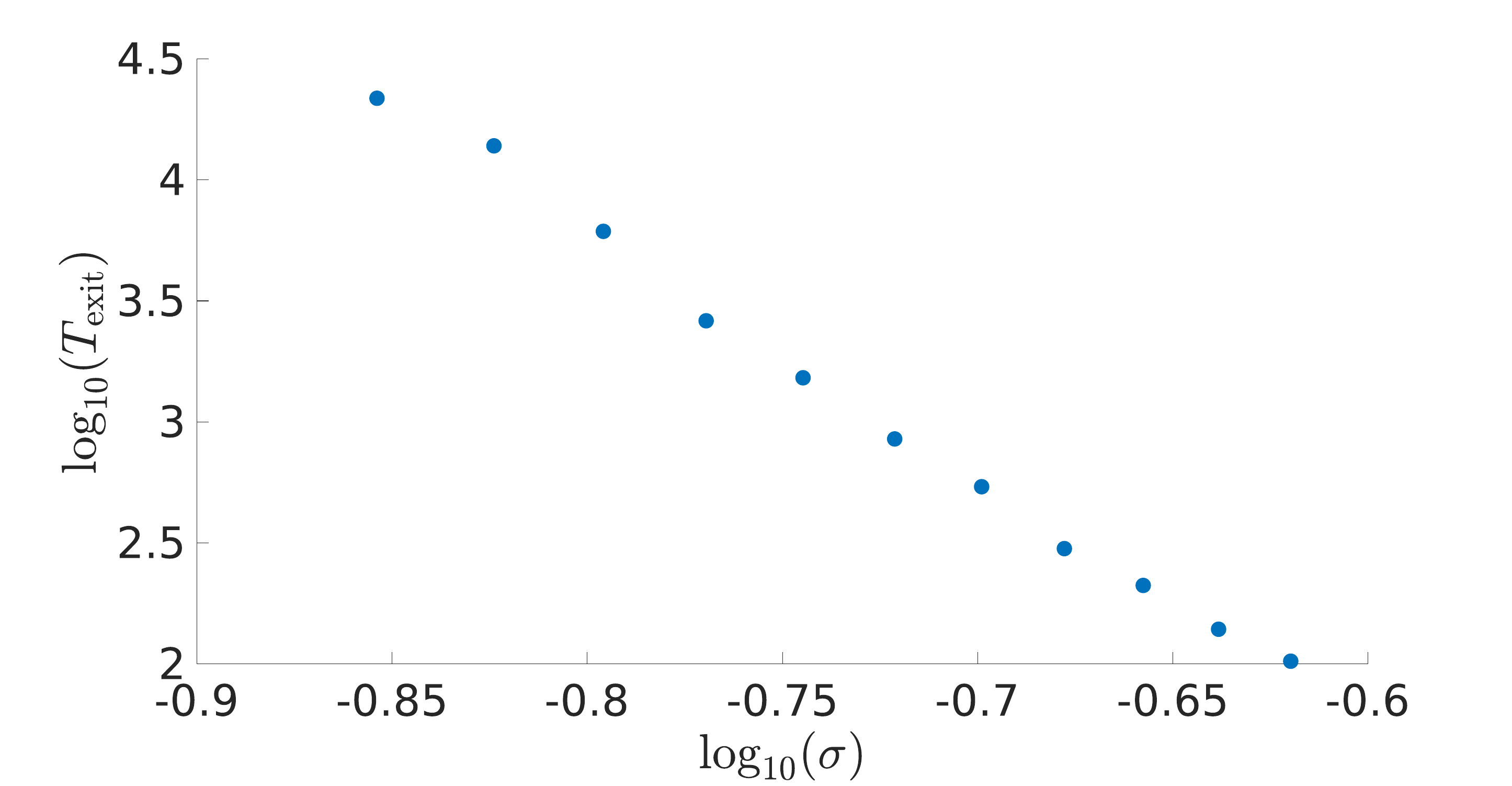
\caption{For a single $(k,a)$-pair where we can properly resolve $E[T_{\rm{exit}}(k,a,\sigma)]$, we can study the relation between $E[T_{\rm{exit}}]$ and $\sigma$. Here, we choose $(k,a)=(30,2)$. The value $k=30$ was deliberately chosen as it is the value for which $E[T_{\rm{exit}}(k,2)]$ attains its maximum for small $\sigma$-values. However, For $\sigma >0.2$, $k=31$ becomes the wave number with the longest average exit time. We plot on a $\log-\log$ scale the average first exit time as a function of $\sigma$ for $a=2$ and $k=30$. This figure indicates that $T_{\rm{exit}}(30,2,\sigma)\sim \sigma^\alpha$ for $\alpha\approx10$. }
    \label{fig:app:AFETvsSigma}
\end{figure}

\end{document}